\documentclass[11pt]{amsart}
\usepackage{amssymb}
\usepackage{epsfig}
\usepackage{mathdots}
\usepackage{hyperref}
\usepackage{fullpage}
\usepackage{color}

\usepackage{array, longtable}
\newcolumntype{C}{>{$}c<{$}}

 \theoremstyle{plain}
\newtheorem{theorem}{Theorem}[section]
\newtheorem{corollary}{Corollary}[section]

\newtheorem{lemma}{Lemma}[section]
\newtheorem{proposition}{Proposition}[section]
\newtheorem{example}{Example}

\theoremstyle{definition}
\newtheorem{definition}{Definition}
\newtheorem{remark}{Remark}[section]

\numberwithin{equation}{section}

\setbox0=\hbox{$+$}
\newdimen\plusheight
\plusheight=\ht0
\def\+{\;\lower\plusheight\hbox{$+$}\;}

\setbox0=\hbox{$-$}
\newdimen\minusheight
\minusheight=\ht0
\def\-{\;\lower\minusheight\hbox{$-$}\;}

\setbox0=\hbox{$\cdots$}
\newdimen\cdotsheight
\cdotsheight=\plusheight
\def\cds{\lower\cdotsheight\hbox{$\cdots$}}

\newcommand{\lra}{\Longrightarrow}

\newcommand{\lrla}{\Longleftrightarrow}

\begin{document}

\title[Identical Vanishing of Coefficients modulo 4, 9 and 25]{Identical Vanishing of Coefficients in the Series Expansion of Eta Quotients, modulo 4, 9 and 25}
\author{Tim Huber, James McLaughlin, and Dongxi Ye}

\address{
School of Mathematical and Statistical Sciences, University of Texas Rio Grande
Valley, Edinburg, Texas 78539, USA}
\email{timothy.huber@utrgv.edu}

\address{Mathematics Department,
 25 University Avenue,
West Chester University, West Chester, PA 19383}
\email{jmclaughlin2@wcupa.edu}

\address{
School of Mathematics (Zhuhai), Sun Yat-sen University, Zhuhai 519082, Guangdong,
People's Republic of China}

\email{yedx3@mail.sysu.edu.cn}

\thanks{Dongxi Ye was supported by the Guangdong Basic and Applied Basic Research
Foundation (Grant No. 2024A1515030222). }

\date{\today}

\begin{abstract}
Let $A(q)=:\sum_{n=0}^{\infty}a_n q^n$ and $B(q)=:\sum_{n=0}^{\infty}b_n q^n$ be two eta quotients. In some previous papers, the
present authors considered the problem of when
\[
a_n=0 \lrla b_n=0.
\]
In the present paper we consider the ``mod $m$'' version of this problem, i.e.  for which eta quotients $A(q)$ and $B(q)$ and for which integers $m>1$ do we have (non-trivially) that
\[
a_n \equiv 0 \pmod m \lrla b_n \equiv 0 \pmod m?
\]
(We say ``non-trivially'' as there are trivial situations where $a_n \equiv b_n \pmod m$ for all $n\geq 0$).

The $m$ for which we found non-trivial (in the sense just mentioned) results were $m=p^2$, $p=2, 3$ and $5$. For $m=4$ and $m=9$, we found  results which apply  to infinite families of eta quotients. One such is the following:
 Let $A(q)$ be any eta quotient of the form
  \begin{equation*}
  A(q) = f_1^{3j_1+1}\prod_{3\nmid i}f_i^{3j_i}\prod_{3|i}f_i^{j_i} =: \sum_{n=0}^{\infty}a_nq^n,
  \hspace{25pt}
   B(q) = \frac{f_3}{f_1^3}A(q) =: \sum_{n=0}^{\infty}b_nq^n
  \end{equation*}
  with $f_{k}=\prod_{n=1}^{\infty}(1-q^{kn})$.
  Then
  \begin{align*}
    a_{3n}- b_{3n}  & \equiv 0\pmod 9, \\
   2a_{3n+1}+b_{3n+1} &\equiv 0 \pmod 9, \\
   a_{3n+2}+2b_{3n+2} &\equiv 0 \pmod 9.
  \end{align*}
 Some of these theorems also had some combinatorial implications, one example being the following:  Let $p_2^{(3)}(n)$ denote the number of bipartitions $(\pi_1, \pi_2)$ of $n$ where $\pi_1$ is 3-regular. Then
  \begin{equation*}
  p_2^{(3)}(n)\equiv 0 \pmod 9 \lrla n \text{ is \underline{not} a generalized pentagonal number}.
\end{equation*}
In the case of $m=25$, we do not have any general theorems that apply to an infinite family of eta quotients, such as the modulo 9 result stated above. Instead we give two tables of results that appear to hold experimentally. Proofs of results stated in these tables appear to need the theory of modular forms and are more complicated. We do prove some individual results, such as the following:
 Let the sequences  $\{c_n\}$  and  $\{d_n\}$ be defined by
  \begin{equation*}
  f_1^{10}=:\sum_{n=0}^{\infty}c_nq^n, \hspace{25pt}  f_1^{5}f_5=:\sum_{n=0}^{\infty}d_nq^n.
  \end{equation*}
  Then
  \begin{equation*}
    c_n \equiv 0 \pmod{25} \lrla d_n \equiv 0 \pmod{25}.
  \end{equation*}
\end{abstract}

\maketitle

\allowdisplaybreaks




\section{Introduction}
As usual, for $|q|<1$, define
\begin{align*}
   f_j &:=     (q^{j};q^{j})_\infty := \prod_{n=1}^{\infty}(1-q^{jn}),
\end{align*}
and recall that an eta quotient is a finite product of the form
\[
\prod_{j}f_j^{n_j},
\]
for some
  integers $j \in \mathbb{N}$ and  $n_j \in \mathbb{Z}$.

For an eta quotient $A(q)$ with series expansion  $A(q)=\sum_{n\geq 0} a_n q^n$, define
\[
A_{(0)}:=\{n\in \mathbb{N}: a_n=0\}.
\]
If $A(q)$ and $B(q)$ are two eta quotients for which $A_{(0)}=B_{(0)}$, then  we say that  $A(q)$ and $B(q)$ have \emph{identically vanishing coefficients}. In some previous papers (\cite{HMcLY22}, \cite{HLMcLYYZ}, \cite{HMcLY23}, \cite{HMcLY23a}) the authors investigated this phenomenon, proving many cases and conjecturing many more.

In the present paper, we investigate the ``modulo $m$'' version of this phenomenon. Let $A(q)=\sum_{n=0}^{\infty}a_nq^n$ and $B(q)=\sum_{n=0}^{\infty}b_nq^n$ be eta quotients, and let $m>1$ be positive integers. We are interested in the situation where
\begin{equation}\label{an0bn0mod}
  a_n \equiv 0 \pmod {m}\lrla b_n \equiv 0 \pmod {m}.
\end{equation}

We first discount some trivial situations where \eqref{an0bn0mod} holds.
If $a_n \equiv b_n \pmod m$ for all $n \geq 0$, then for ease of notation we write
\begin{equation}\label{ABmodmtriv}
A(q)\equiv B(q) \pmod m.
\end{equation}
There are many instances where this holds, and if \eqref{ABmodmtriv} holds, then  \eqref{an0bn0mod} holds trivially. For example,
if $p$ is a prime and
\[
B(q) = \frac{f_1^p}{f_p} A(q) \lra B(q)\equiv A(q) \pmod p, \text{ since }\frac{f_1^p}{f_p}\equiv 1 \pmod p.
\]
A second situation where \eqref{an0bn0mod} may hold trivially occurs if $A(q)$ and $B(q)$ have \emph{similar $m$-dissections}.
\begin{definition}
  By the $m$-dissection of a function $G(q)=\sum_{n=0}^{\infty}g_nq^n$ we mean an expansion of the form
  \begin{equation}\label{mdissecteq1}
   G(q) = \gamma_0G_0(q^m)+\gamma_1 q G_1(q^m)+ \dots + \gamma_{m-1} q^{m-1} G_{m-1}(q^m),
  \end{equation}
  where each dissection component $G_i(q^m)$ is not identically zero ($\gamma_i=0$ is allowed).
  In other words, for each $i$, $0 \leq i \leq m-1$,
  \begin{equation}\label{mdissecteq2}
  \gamma_i q^i G_i(q^m) = \sum_{n=0}^{\infty}g_{mn+i}q^{mn}.
  \end{equation}
\end{definition}
Now suppose $A(q)$ and $B(q)$ are two eta quotients  whose $m$-dissections are given by
\begin{align}\label{ABdissecteq1}
  A(q)& = c_0G_0(q^m)+c_1 q G_1(q^m)+ \dots + c_{m-1} q^{m-1} G_{m-1}(q^m),\\
  B(q)&  = d_0G_0(q^m)+d_1 q G_1(q^m)+ \dots + d_{m-1} q^{m-1} G_{m-1}(q^m), \notag
\end{align}
where $c_i=0\lrla d_i=0$, $i=0,1,\dots, m-1$. It can be seen that if the non-zero $c_i$ and $d_i$ are relatively prime to $m$, then once again \eqref{an0bn0mod} holds trivially.

Thus the situation we are interested in is where we have pairs of eta quotients $A(q)$ and $B(q)$ for which it is the case that they neither have similar $m$-dissections nor are such that \eqref{ABmodmtriv} holds, but \eqref{an0bn0mod} does hold.

For ease of writing and to allow for a slight generalization of the discussion, we introduce some additional notation, and for the eta quotient $A(q)$ as above and $0\leq i <m$, define
\begin{equation}\label{Aimeq}
  A_{i,m}:=\{n\in \mathbb{N}_0|a_n \equiv i \pmod m\}.
\end{equation}
With this notation, \eqref{an0bn0mod} could be rewritten as
\begin{equation}\label{an0bn0mod2}
  A_{0,m}=B_{0,m}.
\end{equation}

The non-trivial cases we found all involve $m$ of the form
 $m=p^2$, specifically  $p=2$, $3$ and $5$.
We next briefly summarize some of the results in the paper.

\subsection{The case \texorpdfstring{$p=2$}{}}
One result proved in the paper is contained in the following theorem.
\begin{theorem}
Let
  \begin{equation*}
   A(q) = \prod f_j^{n_j}=:\sum_{n=0}^{\infty}a_n q^n
  \end{equation*}
  be any eta quotient satisfying the following conditions:
  \begin{itemize}
    \item $n_1$ is odd.
    \item If $j>1$ is odd, then $n_j$ is even.
  \end{itemize}
  Let
 \begin{equation*}
  B(q) = A(q) \frac{f_1^2f_2}{f_4}=:\sum_{n=0}^{\infty}b_n q^n.
 \end{equation*}
 Then for all integers $n\geq 0$,
 \begin{align*}
 a_{2n}-b_{2n} &\equiv 0 \pmod{4} \\
 a_{2n+1}+b_{2n+1} &\equiv 0 \pmod{4}. \notag
 \end{align*}
 \end{theorem}
This theorem also has some combinatorial applications.  For example, we recover the following result of Merca.
\begin{corollary}(Merca \cite[Page 121, Cor. 1]{M21})
 Let $n$ be a positive integer. The number of representations of $n$ as the sum of a generalized
pentagonal number and a square or a twice square is odd if and only if $n$ is an odd generalized pentagonal
number.
\end{corollary}
A second application involves the partition function $p(n)$.
\begin{corollary}\label{exampleforp}
Recall that $p(n)$ denotes the number of unrestricted partitions of the integer $n$, and let
\[
S_{\square}:=\{n^2:n\geq 1\}\cup \{2n^2:n\geq 1\}.
\]
(i)  For any $N$ with $p(N)$ even, or any even $N$ with $p(N)$ odd, one has that
  \begin{equation*}
 \#\{(m,n)| p(m) \text{ odd }, n\in S_{\square}, m+n=N\}
  \end{equation*}
is even.

(ii) For any odd $N$ with $p(N)$ odd  one has that
 \begin{equation*}
   \#\{(m,n)| p(m) \text{ odd }, n\in S_{\square}, m+n=N\}
 \end{equation*}
is odd.
\end{corollary}
{One shall see after Corollary~\ref{ex2} that similar statements hold about the number of $t$-cores and $t$-regular partitions, in the case where $t$ is even, and also about number of partitions into distinct/odd parts.}

We also prove a second general class of congruence results which also has combinatorial applications, one example involving the representation of a positive integer as a sum of a generalized pentagonal number plus three times a square, and a second example involving partitions into distinct parts.

\subsection{The case \texorpdfstring{$p=3$}{}}

Modulo 9, our main result is contained in the following theorem.

\begin{theorem}
  Let $A(q)$ be any eta quotient of the form
  \begin{equation*}
  A(q) = f_1^{3j_1+1}\prod_{3\nmid i}f_i^{3j_i}\prod_{3|i}f_i^{j_i} =: \sum_{n=0}^{\infty}a_nq^n.
  \end{equation*}
  Let
  \begin{equation*}
   B(q) = \frac{f_3}{f_1^3}A(q) =: \sum_{n=0}^{\infty}b_nq^n.
  \end{equation*}
  Then
  \begin{align*}
    a_{3n}- b_{3n}  & \equiv 0\pmod 9, \\
   2a_{3n+1}+b_{3n+1} &\equiv 0 \pmod 9, \\
   a_{3n+2}+2b_{3n+2} &\equiv 0 \pmod 9.
  \end{align*}
\end{theorem}
As an application of this theorem, we have the following result.
\begin{corollary}
  Let $p_2^{(3)}(n)$ denote the number of bipartitions $(\pi_1, \pi_2)$ of $n$ where $\pi_1$ is 3-regular. Then
  \begin{equation*}
  p_2^{(3)}(n)\equiv 0 \pmod 9 \lrla n \text{ is \underline{not} a generalized pentagonal number}.
\end{equation*}
\end{corollary}

\subsection{The case \texorpdfstring{$p=5$}{}}

Modulo 25, we did not find any general infinite families of results, such as in the theorems above. The results in this section are much more experimental, and we list quintuples of eta quotients
\[
\left(\frac{f_1^5}{f_5}\right)^j F(q), \hspace{25pt} 0 \leq j \leq 4.
\]
for which there appears to be results concerning identical vanishing of coefficients modulo 25 (clearly there is identical vanishing of coefficients modulo 5).

We prove two results that were found experimentally. It will be seen that the proofs here are much more technical, and essentially involves extending the methods used in \cite{HMcLY22}. The difficulty of these proofs indicate the desirability of finding more efficient methods of proof. We prove the following two results.
\begin{theorem}
   Let the sequences  $\{c_n\}$  and  $\{d_n\}$ be defined by
  \begin{equation*}
  f_1^{10}=:\sum_{n=0}^{\infty}c_nq^n, \hspace{25pt}  f_1^{5}f_5=:\sum_{n=0}^{\infty}d_nq^n.
  \end{equation*}
  Then
  \begin{equation*}
    c_n \equiv 0 \pmod{25} \lrla d_n \equiv 0 \pmod{25}.
  \end{equation*}
\end{theorem}

\begin{theorem}
    Let $a(n)$ and $b(n)$ be defined by
    \[
    f_1f_5 = \sum_{n=0}^{\infty}a_{n}q^n, \hspace{25pt} f_1^6 =  \sum_{n=0}^{\infty}b_{n}
    q^n.
    \]
    Then
    $$
    \{n|\,b_{n}\equiv0\pmod{25}\}\subsetneqq \{n|\,a_{n}\equiv0\pmod{25}\}.
    $$
\end{theorem}

It is hoped that the experimental data in the tables will provide some insight to others to derive similar results, or possibly additional methods of proving the results which the data appear to suggest.

\section{Vanishing Coefficients Modulo 4
}
In this section we consider pairs of eta quotients $(A(q), B(q))$ for which $A(q) \not \equiv B(q) \pmod 4$, but $A_{0,4}\equiv B_{0,4}$. We prove two general theorems, and also consider some numerical/combinatorial consequences of these.

 Before coming to the first of these, we recall the notation (see, for example, \cite{H88}), for $a$ an integer and $m$ a positive integer ,
\begin{equation}\label{jbardef}
J_{a,m}:=(q^a,q^{m-a},q^m;q^m)_{\infty}, \hspace{25pt}
\bar{J}_{a,m}:=(-q^a,-q^{m-a},q^m;q^m)_{\infty}.
\end{equation}
This notation appears in the  2-dissections of $f_1$ and $1/f_1$.
\begin{lemma}
  The following 2-dissections hold.
  \begin{align}
   f_1&=\frac{f_2 }{f_4}\left(\bar{J}_{6,16}-q \bar{J}_{2,16}\right), \label{f1eqsa}\\
   \frac{1}{f_1}&=\frac{1}{f_2^2}\left(\bar{J}_{6,16}+q \bar{J}_{2,16}\ \right).\label{f1eqsb}
  \end{align}
\end{lemma}

\begin{proof}
  The second identity \eqref{f1eqsb} was proven by  Hirschhorn \cite[Lemma 1]{H160}, and \eqref{f1eqsa} is its $q\to -q$ partner; that is, the reformulation of the identitity upon replacing $q\to -q$ and using the fact
\begin{equation}\label{qmin}
(-q;-q)_{\infty}=\frac{(q^2;q^2)^3_{\infty}}{(q;q)_{\infty}(q^4;q^4)_{\infty}}=\frac{f_2^3}{f_1f_4}.
\end{equation}
\end{proof}

\begin{lemma}\label{jbarcongl}
  The following congruences hold:
  \begin{align}
   -\frac{2 q^2  f_{16}^2 \bar{J}_{2,16}}{ f_8}-\frac{f_2^2 \bar{J}_{6,16}}{f_4}+2
   \bar{J}_{6,16}-\frac{ f_8^5 \bar{J}_{6,16}}{f_4^2 f_{16}^2}&\equiv 0 \pmod 4,\label{jbarcong1}\\
\frac{ f_2^2 \bar{J}_{2,16}}{f_4}+\frac{ f_8^5 \bar{J}_{2,16}}{f_4^2 f_{16}^2}+\frac{2 f_{16}^2
   \bar{J}_{6,16}}{ f_8} &\equiv 0 \pmod 4. \label{jbarcong3}
  \end{align}
\end{lemma}
\begin{proof}
  These can be verified by using modular properties of quotients of Klein forms and Sturm's theorem.

  For~\eqref{jbarcong3}, first notice that
  $$
  \bar{J}_{2,16}=\frac{J_{2,16}(q^{2})}{J_{2,16}(q)}\frac{f_{16}^{2}}{f_{32}}\quad\mbox{and}\quad \bar{J}_{6,16}=\frac{J_{6,16}(q^{2})}{J_{6,16}(q)}\frac{f_{16}^{2}}{f_{32}}.
  $$
  Also, setting $q=e^{2\pi i\tau}$ for ${\rm Im}(\tau)>0$, one can find that
$$
F(\tau)=q^{\frac{7}{8}}\frac{f_{16}^{3}}{J_{2,16}(q)}\quad\mbox{and}\quad G(\tau)=q^{\frac{15}{8}}\frac{f_{16}^{3}}{J_{6,16}(q)}
$$
are both holomorphic modular forms of weight~$1$ for $\Gamma(16)$. See, e.g., \cite[Lemma~2.1]{HMOY}. It is easy to see after some simple manipulations that~\eqref{jbarcong3} amounts to
$$
q^{2}f_{2}^{2}f_{4}f_{8}f_{16}^{2}G(2\tau)F(\tau)+q^{2}f_{8}^{6}G(2\tau)F(\tau)+2q^{3}f_{4}^{2}f_{16}^{4}F(2\tau)G(\tau)\equiv0\pmod{4},
$$
where the left hand side is a holomorphic modular form of weight~$5$ for $\Gamma(32)$. The verification of the congruence can be done by routine computation through an iteration of Sturm's theorem modulo $2$.

The proof of~\eqref{jbarcong1} is similar, but one has to replace $2\bar{J}_{6,16}$ with $2\frac{f_{2}^{5}}{f_{1}^{2}f_{4}^{2}}\bar{J}_{6,16}$ by the fact
$$
2\equiv 2\sum_{n=-\infty}^{\infty}q^{n^{2}}=2\frac{f_{2}^{5}}{f_{1}^{2}f_{4}^{2}} \pmod{4}
$$
in order to construct a modular form of the claimed weight.
\end{proof}

 \begin{theorem}\label{mod4t1}
Let
  \begin{equation}\label{Aqt1eq}
   A(q) = \prod f_j^{n_j}=:\sum_{n=0}^{\infty}a_n q^n
  \end{equation}
  be any eta quotient satisfying the following conditions:
  \begin{itemize}
    \item $n_1$ is odd.
    \item If $j>1$ is odd, then $n_j$ is even.
  \end{itemize}
  Let
 \begin{equation}\label{Bqt1eq}
  B(q) = A(q) \frac{f_1^2f_2}{f_4}=:\sum_{n=0}^{\infty}b_n q^n.
 \end{equation}
 Then for all integers $n\geq 0$,
 \begin{align}\label{anbnt1eq}
 a_{2n}-b_{2n} &\equiv 0 \pmod{4} \\
 a_{2n+1}+b_{2n+1} &\equiv 0 \pmod{4}. \notag
 \end{align}
 \end{theorem}

 \begin{proof}
   Recall that the exponent of $f_1$ in $A(q)$ is odd, so that the exponent of $f_1$ in $A(q)/f_1$ is even, and let the 2-dissection of the latter eta quotient be
   \begin{equation}\label{Aq2dissect}
   \frac{A(q)}{f_1}=: A_0(q^2) + q A_1(q^2) =: A_0 + q A_1.
   \end{equation}
Note for later use that, from the definition of $A(q)$, all the $f_j$ with $j$ odd that occur in $A(q)/f_1$ occur with even exponent, and hence all the coefficients in $qA_1$ are even. Hence, using \eqref{f1eqsa}, the 2-dissection of $A(q)$ is given by
\begin{multline}\label{Aq2disseq}
A(q) = f_1\frac{A(q)}{f_1} = \frac{f_2 }{f_4}\left(\bar{J}_{6,16}-q \bar{J}_{2,16}\right)(A_0 + q A_1)\\
=\left(\frac{A_0 f_2\bar{J}_{6,16}}{f_4}-\frac{q^2 A_1 f_2 \bar{J}_{2,16}}{f_4}\right)+q \left(-\frac{A_0 f_2 \bar{J}_{2,16}}{f_4}+\frac{A_1 f_2 \bar{J}_{6,16}}{f_4}\right).
\end{multline}
Before coming to the 2-dissection of $B(q)$, recall that
\begin{equation}\label{f12f2eq}
  \frac{f_{1}^{2}}{f_{2}}=\sum_{n=-\infty}^{\infty}(-1)^{n}q^{n^{2}}=1+2\sum_{n=1}^{\infty}(-1)^{n}q^{n^{2}},
\end{equation}
a special case of the Jacobi triple product identity. This implies that
\begin{equation}\label{f12f2overf4cong}
 \frac{f_1^2f_2}{f_4}=\frac{f_1^2}{f_2}\frac{f_2^2}{f_4}
=\left(\frac{f_1^2}{f_2}-1\right) \left(\frac{f_2^2}{f_4}-1\right)+\frac{f_1^2}{f_2}+\frac{f_2^2}{f_4}-1
\equiv \frac{f_1^2}{f_2}+\frac{f_2^2}{f_4}-1 \pmod 4.
\end{equation}

Since the 2-dissection of $f_1^2/f_2$ is given by
\[
\frac{f_1^2}{f_2}=\frac{f_8^5}{f_4^2f_{16}^2} -2 q \frac{f_{16}^2}{f_8},
\]
one has that, modulo 4,
\begin{multline}\label{Bq2disseq}
 B(q) \equiv
   \frac{f_2 }{f_4}\left(\bar{J}_{6,16}-q \bar{J}_{2,16}\right)(A_0 + q A_1)
   \left(\frac{f_8^5}{f_4^2f_{16}^2} -2 q \frac{f_{16}^2}{f_8} +\frac{f_2^2}{f_4}-1\right)\\
   =
   \left(\frac{A_0 f_2^3 \bar{J}_{6,16}}{f_4^2}-\frac{A_0 f_2 \bar{J}_{6,16}}{f_4}+\frac{A_0 f_2 f_8^5 \bar{J}_{6,16}}{f_4^3
   f_{16}^2}\right)
   +q \bigg(-\frac{A_0 f_2^3 \bar{J}_{2,16}}{f_4^2}+\frac{A_0 f_2 \bar{J}_{2,16}}{f_4}-\frac{A_0 f_2 f_8^5
   \bar{J}_{2,16}}{f_4^3 f_{16}^2}\\
 +\frac{A_1 f_2^3 \bar{J}_{6,16}}{f_4^2}  -\frac{A_1 f_2 \bar{J}_{6,16}}{f_4}+\frac{A_1 f_2 f_8^5
   \bar{J}_{6,16}}{f_4^3 f_{16}^2}-\frac{2 A_0 f_2 f_{16}^2 \bar{J}_{6,16}}{f_4 f_8}\bigg)\\
   +q^2 \left(-\frac{A_1 f_2^3
   \bar{J}_{2,16}}{f_4^2}+\frac{A_1 f_2 \bar{J}_{2,16}}{f_4}-\frac{A_1 f_2 f_8^5 \bar{J}_{2,16}}{f_4^3 f_{16}^2}+\frac{2 A_0 f_2 f_{16}^2
   \bar{J}_{2,16}}{f_4 f_8}-\frac{2 A_1 f_2 f_{16}^2 \bar{J}_{6,16}}{f_4 f_8}\right)\\
   +\frac{2 q^3 A_1 f_2 f_{16}^2 \bar{J}_{2,16}}{f_4
   f_8}.
\end{multline}
From the dissections at \eqref{Aq2disseq} and \eqref{Bq2disseq}, one  gets that
\begin{multline}\label{a2nb2nmod4}
 \sum_{n=0}^{\infty}(a_{2n}-b_{2n})q^{2n} \equiv
 \left(\frac{A_0 f_2\bar{J}_{6,16}}{f_4}-\frac{q^2 A_1 f_2 \bar{J}_{2,16}}{f_4}\right)\\
  -
 \bigg[
  \left(\frac{A_0 f_2^3 \bar{J}_{6,16}}{f_4^2}-\frac{A_0 f_2 \bar{J}_{6,16}}{f_4}+\frac{A_0 f_2 f_8^5 \bar{J}_{6,16}}{f_4^3
   f_{16}^2}\right)\\
   +q^2 \left(-\frac{A_1 f_2^3
   \bar{J}_{2,16}}{f_4^2}+\frac{A_1 f_2 \bar{J}_{2,16}}{f_4}-\frac{A_1 f_2 f_8^5 \bar{J}_{2,16}}{f_4^3 f_{16}^2}+\frac{2 A_0 f_2 f_{16}^2
   \bar{J}_{2,16}}{f_4 f_8}-\frac{2 A_1 f_2 f_{16}^2 \bar{J}_{6,16}}{f_4 f_8}\right)
 \bigg]\\
 =\frac{f_2 A_0}{f_4} \left(-\frac{2 q^2  f_{16}^2 \bar{J}_{2,16}}{ f_8}-\frac{f_2^2 \bar{J}_{6,16}}{f_4}+2
   \bar{J}_{6,16}-\frac{ f_8^5 \bar{J}_{6,16}}{f_4^2 f_{16}^2}\right)\\
   +\frac{q^2 f_2A_1}{f_4} \left(\frac{ f_2^2 \bar{J}_{2,16}}{f_4}-2
  \bar{J}_{2,16}+\frac{ f_8^5 \bar{J}_{2,16}}{f_4^2 f_{16}^2}+\frac{2  f_{16}^2 \bar{J}_{6,16}}{ f_8}\right)
   \equiv 0 \pmod4,
\end{multline}
where the last congruences follow from \eqref{jbarcong1} for the factor multiplying $f_2A_0/f_4$. For the factor multiplying $q^2 f_2 A_1/f_4$, as mentioned above the coefficients of $A_1$ are all even, and thus all that is necessary is to show that
 \[
 \frac{ f_2^2}{f_4}+\frac{f_8^5 }{f_4^2 f_{16}^2}\equiv 0 \pmod 2,
 \]
 which follows from \eqref{f12f2eq} and the expansion
 \[
 \frac{f_{2}^{5}}{f_{1}^{2}f_{4}^{2}}=\sum_{n=-\infty}^{\infty}q^{n^{2}}=1+2\sum_{n=1}^{\infty}q^{n^{2}}.
 \]

Likewise from the dissections at \eqref{Aq2disseq} and \eqref{Bq2disseq}, one  has
\begin{multline}\label{a2n1b2n1mod4}
 \sum_{n=0}^{\infty}(a_{2n+1}+b_{2n+1})q^{2n+1}\\
  \equiv
 q \left(-\frac{A_0 f_2 \bar{J}_{2,16}}{f_4}+\frac{A_1 f_2 \bar{J}_{6,16}}{f_4}\right)+
 q \bigg(-\frac{A_0 f_2^3 \bar{J}_{2,16}}{f_4^2}+\frac{A_0 f_2 \bar{J}_{2,16}}{f_4}
 -\frac{A_0 f_2 f_8^5
   \bar{J}_{2,16}}{f_4^3 f_{16}^2}\\
 +\frac{A_1 f_2^3 \bar{J}_{6,16}}{f_4^2}  -\frac{A_1 f_2 \bar{J}_{6,16}}{f_4}+\frac{A_1 f_2 f_8^5
   \bar{J}_{6,16}}{f_4^3 f_{16}^2}-\frac{2 A_0 f_2 f_{16}^2 \bar{J}_{6,16}}{f_4 f_8}\bigg)
   +\frac{2 q^3 A_1 f_2 f_{16}^2 \bar{J}_{2,16}}{f_4
   f_8}\\
   =-\frac{qf_2A_0}{f_4} \left(\frac{ f_2^2 \bar{J}_{2,16}}{f_4}+\frac{ f_8^5 \bar{J}_{2,16}}{f_4^2 f_{16}^2}+\frac{2 f_{16}^2
   \bar{J}_{6,16}}{ f_8}\right)\\
   +\frac{q f_2 A_1}{f_4} \left(\frac{2 q^2 f_{16}^2 \bar{J}_{2,16}}{ f_8}+\frac{ f_2^2
   \bar{J}_{6,16}}{f_4}+\frac{f_8^5 \bar{J}_{6,16}}{f_4^2 f_{16}^2}\right)
   \equiv 0 \pmod 4.
 \end{multline}
 The last congruence uses \eqref{jbarcong3} for the factor multiplying $qf_2A_0/f_4$, and for the factor multiplying $qf_2A_1/f_4$, the argument is identical to that used for the factor multiplying $q^2f_2A_1/f_4$ in the previous paragraph.
 \end{proof}


\begin{corollary}\label{cor1fromt1}
  Let the eta quotients $A(q)$ and $B(q)$,  and the sequences $\{a_n\}$ and $\{b_n\}$ be as in Theorem \ref{mod4t1}.
  Then
  \begin{align}\label{an0bn0mod4c}
    a_n \equiv 0 \pmod 4 &\lrla b_n \equiv 0 \pmod 4, \\
    a_n \equiv 2 \pmod 4 &\lrla b_n \equiv 2 \pmod 4.
  \end{align}
\end{corollary}
\begin{proof}
  This is immediate from \eqref{anbnt1eq}.
\end{proof}

\begin{remark}
    It can be seen that the eta quotients $A(q)$ and $B(q)$ have identically vanishing coefficients modulo~4 non-trivially, since $f_{1}^{2}f_{2}/f_{4}\not\equiv1\pmod4$.
\end{remark}

\begin{corollary}\label{cor2fromt1}
Let $A(q)=\sum_{n\geq 0}a_n q^n$ and  $B(q)=\sum_{n\geq 0}b_n q^n$  be as in Theorem \ref{mod4t1}. Let
 \begin{align*}
    S_o&:=\{n|a_n \equiv 1 \pmod 2\}, \\
    S_e&:=\{n|a_n \equiv 0 \pmod 2 \},\\
    S_{\square}&:=\{n^2:n\geq 1\}\cup \{2n^2:n\geq 1\}.
  \end{align*}
(i)  For any $N \in S_e$, or any even $N$ is $S_o$, one has that
  \begin{equation}\label{cor2fromt1eq1}
 \#\{(m,n)| m\in S_o, n\in S_{\square}, m+n=N\}
  \end{equation}
is even.

(ii) For any odd $N$ in $S_0$  one has that
 \begin{equation}\label{cor2fromt1eq2}
   \#\{(m,n)| m\in S_o, n\in S_{\square}, m+n=N\}
 \end{equation}
is odd.

\end{corollary}

\begin{proof}
From \eqref{Bqt1eq} and  \eqref{f12f2overf4cong},
\begin{align*}
  B(q) &= A(q) \frac{f_1^2f_2}{f_4}\\
  & \equiv A(q) \left[ 1+  \left(\frac{f_1^2}{f_2} -1\right) +\left(\frac{f_2^2}{f_4}-1\right)  \right]\pmod 4\\
  &= A(q) + A(q) \left[2\sum_{n=1}^{\infty}(-1)^n q^{n^2} + 2\sum_{n=1}^{\infty}(-1)^n q^{2n^2}   \right]\\
  &= A(q) + \left[ \sum_{m\in S_o} a_m q^m +  \sum_{m\in S_e} a_m q^m \right]\left[2\sum_{n=1}^{\infty}(-1)^n q^{n^2} + 2\sum_{n=1}^{\infty}(-1)^n q^{2n^2}   \right],\\
  & \equiv A(q) + \left[ \sum_{m\in S_o} a_m q^m  \right]\left[2\sum_{n=1}^{\infty}(-1)^n q^{n^2} + 2\sum_{n=1}^{\infty}(-1)^n q^{2n^2}   \right]\pmod 4,
\end{align*}
where the last congruence follows since if $m\in S_e$, then $2|a_m$.

Now consider a term $a_Nq^N$ in the series expansion of $A(q)$. Since $m\in S_o$, $a_m$ is odd, so that the coefficient $2a_m$ in any term $2a_m q^{m+n^2}$ with $m+n^2=N$, or any term $2a_m q^{m+2n^2}$ with $m+2n^2=N$, is an odd multiple of 2.

It is easy to see from \eqref{anbnt1eq} that for any $N$ satisfying the conditions in part (i) that $b_N \equiv a_N \pmod 4$,
so that there must be an even number of terms of the forms $2a_m q^{m+n^2}$ with $m+n^2=N$ or  $2a_m q^{m+2n^2}$ with $m+2n^2=N$, thus leading to the statement at \eqref{cor2fromt1eq1}.

Likewise from \eqref{anbnt1eq} one has that for any $N$ satisfying the conditions in part (ii) that $b_N \equiv a_N+2 \pmod 4$, so that an odd number of such terms are required, thus giving a proof of the statement at \eqref{cor2fromt1eq2}.
\end{proof}

The previous corollary has some combinatorial implications. The first example is a result of Merca.
\begin{corollary}\label{ex1}(Merca \cite[Page 121, Cor. 1]{M21})
 Let $n$ be a positive integer. The number of representations of $n$ as the sum of a generalized
pentagonal number and a square or a twice square is odd if and only if $n$ is an odd generalized pentagonal
number.
\end{corollary}
\begin{proof}
This follows from Corollary \ref{cor2fromt1}, since $A(q) = f_1$ satisfies the requirements of Theorem~\ref{mod4t1}.
\end{proof}
Our second example involves the partition function $p(n)$.
\begin{corollary}\label{ex2}
Recall that $p(n)$ denotes the number of unrestricted partitions of the integer $n$.

(i)  For any $N$ with $p(N)$ even, or any even $N$ with $p(N)$ odd, one has that
  \begin{equation}\label{cor2fromt1eq1pn}
 \#\{(m,n)| p(m) \text{ odd }, n\in S_{\square}, m+n=N\}
  \end{equation}
is even.

(ii) For any odd $N$ with $p(N)$ odd  one has that
 \begin{equation}\label{cor2fromt1eq2pn}
   \#\{(m,n)| p(m) \text{ odd }, n\in S_{\square}, m+n=N\}
 \end{equation}
is odd.
\end{corollary}
\begin{proof}
This likewise follows  from Corollary \ref{cor2fromt1}, since the generating function for the sequence $p(n)$, namely $A(q)=1/f_1$, satisfies the requirements of Theorem \ref{mod4t1}.
\end{proof}

\begin{example}
(i)  Consider $N=55$, with $p(55)= 451276$, even. The list of pairs $(m,n)$ with $n \in S_{\square}$ such that $m+n=55$ is precisely
  \[
  \{(54, 1 ),  (53, 2 ),  (51, 4 ),  (47, 8 ),  (46, 9 ),  (39, 16 ),  (37, 18 ),  (30, 25 ),  (23, 32 ),  (19, 36 ),  (6, 49 ),  (5, 50 )  \}
  \]
  Since
  \begin{multline*}
    (p(54), p(53), p(51), p(47), p(46), p(39), p(37), p(30), p(23), p(19), p(6), p(5)) \\
   =(386155, 329931, 239943, 124754, 105558, 31185, 21637, 5604, 1255,
490, 11, 7),
  \end{multline*}
  and we retain only those $m$ for which $p(m)$ is odd, the list of pairs $(m,n)$ that satisfy \eqref{cor2fromt1eq1pn} is
  \[
  \{  (54,1), (53,2), (51,4), (39,16), (37,18), (23,32), (6,49), (5,50)\},
  \]
  which has an even number of pairs, namely eight.

  In a similar manner, if one takes $N=60$ (even) so that $p(60)=966467$ (odd), and then proceeds similarly one gets that the list of pairs $(m,n)$ that satisfy \eqref{cor2fromt1eq1pn} is
  \[
  \{  (56,4),(52,8),(51,9),(44,16),(35,25),(24,36)\},
  \]
  which again has an even number of pairs, namely six.

  (ii) On the other hand, if one takes $N=53$ (odd) so that $p(53)=329931$ (odd), and then proceeds similarly one gets that the list of pairs $(m,n)$ that satisfy \eqref{cor2fromt1eq2pn} is
  \[
  \{  (52,1),(51,2),(49,4),(44,9),(37,16),(35,18),(17,36),(4,49),(3,50)\},
  \]
  which  has an odd number of pairs, namely nine.
\end{example}

Note that similar statements to those at \eqref{cor2fromt1eq1pn} and \eqref{cor2fromt1eq2pn} also
hold for $t$-cores  and $t$-regular partitions setting $A(q)=f_{t}^{t}/f_{1}$ and $f_t/f_1$ in Corollary~\ref{cor2fromt1}, respectively, both in the case where $t$ is even, and also for partitions into distinct/odd parts ($f_2/f_1$), since the corresponding generating functions also clearly satisfy the  requirements of Theorem \ref{mod4t1}.

\subsection{Another Infinite Family of Congruences.}

To prove the next family of congruence results, we use the results in the following lemmas.

\begin{lemma}
  The following identity holds.
  \begin{equation}\label{2ndmod4congeq2}
    f_1=\bar{J}_{5,12} - q \bar{J}_{1,12}.
  \end{equation}
\end{lemma}
\begin{proof}
  This follows upon splitting the series representation for $f_1$,
 \begin{equation}\label{f1ser}
   f_1 = \sum_{t=-\infty}^{\infty}(-1)^t q^{t(3t-1)/2}
  \end{equation}
into two series, one with $t$ even and the other with $t$ odd, and then using the Jacobi triple product identity on each of the two series.
\end{proof}
Note for later use that \eqref{2ndmod4congeq2} gives us that
 \begin{equation}\label{2ndmod4congeq3}
     \sum_{\stackrel{t=-\infty}{t \text{ odd}}}^{\infty}(-1)^t q^{t(3t-1)/2}= - q \bar{J}_{1,12}.
  \end{equation}

\begin{lemma}
The following congruence holds.
\begin{equation}\label{2ndmod4congeq1}
 f_1 \left(\frac{f_1^2}{f_2}+\frac{f_3^2}{f_6}\right)-2 \bar{J}_{5,12}\equiv 0 \pmod 4.
\end{equation}
\end{lemma}
\begin{proof}
  The proof is similar to that of Lemma~\ref{jbarcongl}. Noticing that
  $$
  2\equiv 2\sum_{n=-\infty}^{\infty}q^{n^{2}}=2\frac{f_{2}^{5}}{f_{1}^{2}f_{4}^{2}}\pmod{4},
  $$
  $$
  \bar{J}_{5,12}(q)=\frac{J_{5,12}(q^{2})}{J_{5,12}(q)}\frac{f_{12}^{2}}{f_{24}},
  $$
  and
  $$
  F(\tau)=q^{\frac{35}{24}}\frac{f_{12}^{3}}{J_{5,12}(q)}
  $$
  with $q=e^{2\pi i\tau}$ for ${\rm Im}(\tau)>0$ is a holomorphic modular form of weight~$1$ for $\Gamma(12)$, one can see that \eqref{2ndmod4congeq3} is equivalent to
  $$
  q^{\frac{31}{24}}f_{1}^{5}f_{4}^{2}f_{6}f_{12}F(2\tau)+q^{\frac{31}{24}}f_{1}^{3}f_{2}f_{3}^{2}f_{4}^{2}f_{12}F(2\tau)-2q^{\frac{11}{4}}f_{2}^{6}f_{6}f_{24}^{2}F(\tau)\equiv0\pmod{4},
  $$
  where the left hand side is a holomorphic modular form of weight~$\frac{11}{2}$ for $\Gamma(24)$. Finally, square both sides and apply Sturm's theorem to affirm the desired congruence.
\end{proof}

\begin{theorem}\label{conjmod2nn}
  Let $S$ denote any finite set of positive integers and let
  \begin{equation}\label{conjmodnneq1}
    A(q) = f_1 \prod_{j\in S} \left( \frac{f_j^2}{f_{2j}} \right)^{n_j} =: \sum_{n=0}^{\infty}a_nq^n.
  \end{equation}
  Let
  \begin{equation}\label{conjmodnneq11}
  B(q) = A(q) \frac{f_1^2 f_3^2}{f_2f_6} =: \sum_{n=0}^{\infty}b_nq^n.
 \end{equation}
  Then with the notation of \eqref{Aimeq},
  \begin{equation}\label{conjmodnneq2}
   b_n \equiv
   \begin{cases}
   a_n \pmod 4, & a_n \equiv 0 \pmod 2 \text{ or }n =\frac{t(3t-1)}{2},\,\, t \text{ even},\\
    a_n +2 \pmod 4,  & n =\frac{t(3t-1)}{2}, \,\,t \text{ odd}.
    \end{cases}
  \end{equation}
\end{theorem}

\begin{remark}
   Since \eqref{AqBqnew} shows that $A(q) \equiv f_1 \pmod 2$, the expansion at \eqref{f1ser} shows that the cases listed in
\eqref{conjmodnneq2} cover all the possibilities.
Note for what follows that for any integer $j\geq 1$ one has
\begin{equation}\label{fj2overf2jcong}
 \frac{f_j^2}{f_{2j}} \equiv 1 \pmod 2.
\end{equation}
\end{remark}

\begin{proof}
First, once can restrict to the case where each  $n_j=1$ since for any integer $j\geq 1$ one has
\[
\left( \frac{f_j^2}{f_{2j}} \right)^{2}\equiv 1 \pmod 4 \text{ and } \frac{f_{2j}}{f_j^2} \equiv \frac{f_j^2}{f_{2j}} \pmod 4.
\]

Secondly, by \eqref{fj2overf2jcong}, one has that
\begin{align}\label{AqBqnew}
  A(q)& \equiv f_1\left[1 + \sum_{j\in S}\left( \frac{f_j^2}{f_{2j}}-1 \right) \right]\pmod 4,\\
  B(q) & \equiv f_1\left[ 1 + \sum_{j\in S}\left( \frac{f_j^2}{f_{2j}}-1 \right) + \left(\frac{f_1^2}{f_2}-1 \right)
  +\left(\frac{f_3^2}{f_6}-1 \right)\right] \pmod 4\notag\\
  &\equiv A(q) +f_1\left[\left(\frac{f_1^2}{f_2}-1 \right)
  +\left(\frac{f_3^2}{f_6}-1 \right)\right]\pmod 4. \notag
\end{align}
Next, consider the second product,
\begin{align}\label{2ndmod4cneq}
\sum_{n=0}^{\infty}c_nq^n&:= f_1\left[\left(\frac{f_1^2}{f_2}-1 \right)
  +\left(\frac{f_3^2}{f_6}-1 \right)\right]\\
&=f_1 \left(\frac{f_1^2}{f_2}+\frac{f_3^2}{f_6}\right)-2 \bar{J}_{5,12}+ 2q \bar{J}_{1,12} \text{ (by \eqref{2ndmod4congeq2})}\notag\\
& \equiv   2q \bar{J}_{1,12} \pmod 4 \text{ (by \eqref{2ndmod4congeq1})}.\notag
\end{align}
Thus $c_n \equiv 2 \pmod 4$ when $n$ has the form $n=t(3t-1)/2$ with $t$ odd (by \eqref{2ndmod4congeq3}), and $c_n \equiv 0 \pmod 4$ otherwise.

 Since $b_n \equiv a_n + c_n \pmod 4$ by \eqref{AqBqnew}, the result now follows.
\end{proof}

What we have shown implies the following corollary, which parallels the result of Merca in Corollary \ref{ex1}

\begin{corollary}\label{cex1}
  Let $n$ be a positive integer. The number of representations of $n$ as the sum of a generalized
pentagonal number and a square or as the sum of a generalized
pentagonal number and  three times a square is odd if and only if $n$ is a generalized pentagonal
number of the form
\[
n =\frac{t(3t-1)}{2}, \,\, t \text{ odd}.
\]
\end{corollary}
\begin{proof}
  This follows from the series representation of \eqref{2ndmod4cneq} that can be written equivalently as:
  \begin{align}\label{2ndmod4cneq2}
 f_1\left[\left(\frac{f_1^2}{f_2}-1 \right)
  +\left(\frac{f_3^2}{f_6}-1 \right)\right]&\equiv   2q \bar{J}_{1,12} \pmod 4,\\
  \left[ \sum_{m-\infty}^{\infty} (-1)^m q^{m(3m-1)/2}  \right]\left[2\sum_{n=1}^{\infty}(-1)^n q^{n^2} + 2\sum_{n=1}^{\infty}(-1)^n q^{3n^2}   \right]&\equiv   2q \bar{J}_{1,12} \pmod 4 \notag, \\
  \left[ \sum_{m-\infty}^{\infty} (-1)^m q^{m(3m-1)/2}  \right]\left[\sum_{n=1}^{\infty}(-1)^n q^{n^2} + \sum_{n=1}^{\infty}(-1)^n q^{3n^2}   \right]&\equiv   q \bar{J}_{1,12} \pmod 2 \notag,\\
   \left[ \sum_{m-\infty}^{\infty}  q^{m(3m-1)/2}  \right]\left[\sum_{n=1}^{\infty} q^{n^2} + \sum_{n=1}^{\infty} q^{3n^2}   \right]&\equiv   q \bar{J}_{1,12} \pmod 2, \notag
\end{align}
with \eqref{2ndmod4congeq3} being used at the end.
\end{proof}

Theorem \ref{conjmod2nn} also implies the following congruence result for the number of partitions into distinct/odd parts, $Q(n)$.

\begin{corollary}\label{oddpartscong}
 Let $n$ be a positive integer and let $Q(n)$ denote the number of partitions of $n$ into distinct parts. Define
 \[
 v_n =
 \begin{cases}
   1, & \mbox{if } n=\frac{t(3t-1)}{2}, \text{ some }t\in \mathbb{Z}, \\
  0, & \mbox{otherwise},
 \end{cases}
 \]
 and let $w_n$ denote the number of representations of $n$ as a sum of a generalized pentagonal number plus three times a square. Then for any positive integer $N$ one has
 \begin{equation}\label{qncongeq1}
   Q_N \equiv (v_N +2w_N) \pmod 4.
 \end{equation}
\end{corollary}
\begin{proof}
Let
\[
A(q) = f_1\left( \frac{f_1^2}{f_2}\right)^{-1}=\frac{f_2}{f_1}=\sum_{n=0}^{\infty}Q(n)q^n,
\]
with this $A(q)$ satisfying the requirements of Theorem \ref{conjmod2nn}. Then
\begin{align}\label{BqQneq}
  B(q) & =  A(q) \frac{f_1^2 f_3^2}{f_2f_6} = f_1  \frac{f_3^2}{f_6} =  f_1 \left(1 + \frac{f_3^2}{f_6} -1\right)\\
   & =  \sum_{m-\infty}^{\infty} (-1)^m q^{m(3m-1)/2}+\left[ \sum_{m-\infty}^{\infty} (-1)^m q^{m(3m-1)/2}  \right]\left[2 \sum_{n=1}^{\infty}(-1)^n q^{3n^2}   \right]\notag\\
   &\equiv \sum_{m-\infty}^{\infty} (-1)^m q^{m(3m-1)/2}+2\left[ \sum_{m-\infty}^{\infty}  q^{m(3m-1)/2}  \right]\left[ \sum_{n=1}^{\infty} q^{3n^2}   \right]  \pmod 4. \notag
\end{align}
The only subtlety involves the case when $N = t(3t-1)/2$ for some odd integer $t$, in which case the first series contributes a $-1$ to $b_N$ (with the notation of the theorem) and the theorem gives that
\[
a_N = Q(N)\equiv b_N + 2 =-1 + 2w_N +2 = 1+2w_N = v_N + 2w_N \pmod 4,
\]
as in the other cases.
\end{proof}

 \begin{remark}
    This result could be thought of as a refinement of an implication of the Pentagonal Number Theorem, namely that $Q(n)$ is odd if and only if $n$ is a generalized pentagonal number.
 \end{remark}

\section{Vanishing Coefficients Modulo 9 }\label{mod9}

In this section we prove some congruence results modulo 9 for an infinite family of eta quotients. We first prove/state some necessary preliminary results.

\begin{lemma}
The following congruences hold.
    \begin{equation}\label{JJ1}
  \left[2 \left(\frac{f_1^3}{f_3}-1\right) -3q\frac{f_9^3}{f_3}\right]J_{12,27}- 3 q^3 \frac{f_9^3}{f_3}J_{3,27}\equiv 0 \pmod 9.
\end{equation}

\begin{equation}\label{JJ2}
\left[2 \left(\frac{f_1^3}{f_3}-1\right)-3q\frac{f_9^3}{f_3} +3\right]J_{6,27}- 3  \frac{f_9^3}{f_3}J_{12,27}\equiv 0 \pmod 9.
\end{equation}

\begin{equation}\label{JJ3}
\left[2 \left(\frac{f_1^3}{f_3}-1\right)-3q\frac{f_9^3}{f_3} -3\right]J_{3,27}+ 3  \frac{f_9^3}{f_3}J_{6,27}\equiv 0 \pmod 9.
\end{equation}
\end{lemma}

\begin{proof}
    These can be verified by validating the congruences up to Sturm's bound \cite{S87}. We illustrate with~\eqref{JJ1}. Replacing $q$ with $q^{3}$ and multiplying both sides of~\eqref{JJ1} by $q^{4}\frac{f_{81}^{2}}{J_{9,81}J_{36,81}}$, it is equivalent to proving that
    $$
     \left[2 \left(\frac{f_3^3}{f_9}-1\right) -3q^{3}\frac{f_{27}^3}{f_9}\right]q^{4}\frac{f_{81}^{2}}{J_{9,81}}- 3 q^{3} \frac{f_{27}^3}{f_9}\left(q^{10}\frac{f_{81}^{2}}{J_{36,81}}\right)\equiv 0 \pmod 9.
    $$
    In terms of $\tau$ via $q=\exp(2\pi i\tau)$, each of
    $$
    \frac{f_3^3}{f_9},\, q^{3}\frac{f_{27}^3}{f_9},\, q^{4}\frac{f_{81}^{2}}{J_{9,81}},\, q^{10}\frac{f_{81}^{2}}{J_{36,81}}
    $$
    can be seen to be a holomorphic modular form of weight~$1$ for $\Gamma_{1}(81)$ using well known facts about the modularity of Klein forms and eta quotients. Therefore, the form on the left hand side of the congruence is a holomorphic modular form of weight~$2$ for $\Gamma_{1}(81)$. A double iteration of Sturm's theorem \cite{S87} modulo $3$ asserts that for $a(n)$ defined by
    $$
    \sum_{n=0}^{\infty}a(n)q^{n}=\left[2 \left(\frac{f_3^3}{f_9}-1\right) -3q^{3}\frac{f_{27}^3}{f_9}\right]q^{4}\frac{f_{81}^{2}}{J_{9,81}}- 3 q^{3} \frac{f_{27}^3}{f_9}\left(q^{10}\frac{f_{81}^{2}}{J_{36,81}}\right),
    $$
    if $a(n)\equiv0\pmod{9}$ for $n\leq \frac{1}{12}[{\rm SL}_{2}(\mathbb{Z}):\pm\Gamma_{1}(81)]=243$, then $a(n)\equiv0\pmod{9}$ for any $n$. Computationally, one can find that the assumption indeed holds, and thus, the desired congruence follows.

    The other two congruences respectively follow from the following intermediate results
$$
\left[2 \left(\frac{f_3^3}{f_9}-1\right)-3q^{3}\frac{f_{27}^3}{f_9} +3\right]q^{10}\frac{f_{81}^{2}}{J_{36,81}}- 3  q^{3}\frac{f_{27}^3}{f_9}\left(q^{7}\frac{f_{81}^{2}}{J_{18,81}}\right)\equiv0\pmod{9}
$$
and
$$
\left[2 \left(\frac{f_3^3}{f_9}-1\right)-3q^{3}\frac{f_{27}^3}{f_9} -3\right]q^{7}\frac{f_{81}^{2}}{J_{18,81}}+ 3  q^{3}\frac{f_{27}^3}{f_9}\left(q^{4}\frac{f_{81}^{2}}{J_{9,81}}\right)\equiv0\pmod{9},
$$
which can be verified in the same way as the first case.
\end{proof}

Also required is the 3-dissection of $f_1$, which is easily  derived from the Jacobi triple product identity (see also \cite[page 15]{H17}):
\begin{equation}\label{f13diss}
f_1 = J_{12,27}-qJ_{6,27}-q^2J_{3,27},
\end{equation}
where $J_{a,m}$ is as defined at \eqref{jbardef}.

Since $f_3/f_1^3$ appears in Theorem \ref{conjmod91}, we derive a 3-dissection modulo 9 for it.
Recall first that the Borwein theta function $a(q)$ is defined by
\[
a(q)  = \sum_{m,n=-\infty}^{\infty}q^{m^2+mn+n^2}
=1 +6 \sum_{n=1}^{\infty}
\left(
\frac{q^{3n-2}}{1-q^{3n-2}} - \frac{q^{3n-1}}{1-q^{3n-1}}
\right).
\]
It is clear that
\begin{equation}\label{aqm1eq}
 a(q)-1\equiv 0 \pmod 3, \lra a(q)^2 \equiv 2a(q)-1 \pmod 9 \text{ (upon squaring)}.
\end{equation}
 It is also known (see for example, \cite[Eq. (21.3.1), page 183]{H17}) that
\begin{equation}\label{conjmod91eq33}
  a(q^3) = \frac{f_1^3}{f_3}+3q \frac{f_9^3}{f_3}.
\end{equation}

\begin{lemma}
  There holds
  \begin{align}\label{f3f13lem}
   \frac{f_3}{f_1^3}  &\equiv 1 + 3 C_0(q^3)+3 q C_1(q^3) =: 1 + 3 C_0 +3 q C_1 \pmod 9, \text{ where } \\
    3C_0 &\equiv 2 (a(q^3)-1)\pmod 9, \hspace{50pt} C_1 \equiv \frac{f_9^3}{f_3} \pmod 9.
  \end{align}
\end{lemma}

\begin{proof}
By \cite[page 467]{H16}
\begin{align}\label{conjmod91eq4}
\frac{f_3}{f_1^3} &= \frac{f_9^3}{f_3^9}
\left(
a(q^3)^2+3q a(q^3) \frac{f_9^3}{f_3} + 9 q^2 \frac{f_9^6}{f_3^2}
\right)\\
&\equiv a(q^3)^2+3q a(q^3) \frac{f_9^3}{f_3} \pmod 9 \,(\text{since }\frac{f_9^3}{f_3^9} \equiv 1 \pmod 9)\notag\\
& \equiv 2 a(q^3) - 1 +3 q  \frac{f_9^3}{f_3} \pmod 9 \,(\text{by both parts of \eqref{aqm1eq}})\notag\\
&\equiv 2 \frac{f_1^3}{f_3}-1 \pmod 9\, (\text{by \eqref{conjmod91eq33}}) \notag
\end{align}
Thus comparing the last two lines we see that modulo 9,
\begin{align}\label{3C03C1eq1}
  \frac{f_3}{f_1^3}& \equiv 2 \frac{f_1^3}{f_3}-1 =: 1 + 3 C_0(q^3)+3 q C_1(q^3)\pmod 9 =: 1 + 3 C_0 +3 q C_1, \text{ where }\notag\\
3C_0 &\equiv 2 (a(q^3)-1)\pmod 9, \hspace{50pt} C_1 \equiv \frac{f_9^3}{f_3} \pmod 9.
\end{align}

\end{proof}
Note for later use that \eqref{3C03C1eq1} and \eqref{conjmod91eq4} give that
\begin{equation}\label{3C03C1eq2}
 3C_0 \equiv 2 (a(q^3)-1)\equiv 2 \left(\frac{f_1^3}{f_3}-1\right) -3q\frac{f_9^3}{f_3} \pmod 9.
\end{equation}

\begin{theorem}\label{conjmod91}
  Let $A(q)$ be any eta quotient of the form
  \begin{equation}\label{conjmod91eq1}
  A(q) = f_1^{3j_1+1}\prod_{3\nmid i}f_i^{3j_i}\prod_{3|i}f_i^{j_i} =: \sum_{n=0}^{\infty}a_nq^n.
  \end{equation}
  Let
  \begin{equation}\label{conjmod91eq2}
   B(q) = \frac{f_3}{f_1^3}A(q) =: \sum_{n=0}^{\infty}b_nq^n.
  \end{equation}
  Then
  \begin{align}
    a_{3n}- b_{3n}  & \equiv 0\pmod 9, \label{conjmodn1eq3a}\\
   2a_{3n+1}+b_{3n+1} &\equiv 0 \pmod 9, \label{conjmodn1eq3b}\\
   a_{3n+2}+2b_{3n+2} &\equiv 0 \pmod 9.\label{conjmodn1eq3c}
  \end{align}
\end{theorem}

\begin{proof}
For ease of notation we rewrite the $q$-products at \eqref{conjmod91eq1} by defining the functions $A_0$, $A_1$, $A_2$ and $D$ from their 3-dissections by setting
\begin{align*}
  f_1^{3j_1}\prod_{3\nmid i}f_i^{3j_i} & =: A_0(q^3)+3q A_1(q^3) + 3 q^2 A_2(q^3)=: A_0+3qA_1+3q^2 A_2, \\
  \prod_{3|i}f_i^{j_i} & =: D(q^3)=:D.
\end{align*}
Then
\begin{multline}\label{A(q)}
 A(q) = (J_{12,27}-qJ_{6,27}-q^2J_{3,27}) (A_0+3qA_1+3q^2 A_2)D\\
= D\left[A_0 J_{12,27}- 3q^3 \left( A_1  J_{3,27}+ A_2 J_{6,27}\right) \right]\\
  +q D\left[ 3 A_1  J_{12,27}-A_0 J_{6,27}  -3 q^3 A_2   J_{3,27}\right]
 +q^2 D \left[-A_0  J_{3,27}-3 A_1   J_{6,27}+3 A_2  J_{12,27}\right].
\end{multline}
Similarly, one has that
\begin{multline}\label{B(q)}
 B(q) \equiv  (J_{12,27}-qJ_{6,27}-q^2J_{3,27}) (A_0+3qA_1+3q^2 A_2)D(1 + 3 C_0 +3 q C_1)\\
= D(A_0 J_{12,27}+3 A_0 C_0 J_{12,27})\\
+q D \left(-A_0
   J_{6,27}-3 A_0 C_0 J_{6,27}+3 A_1 J_{12,27}+9
   A_1 C_0 J_{12,27}+3 A_0 C_1
   J_{12,27}\right)\\
   +q^2 D \big(-A_0 J_{3,27}-3 A_0
   C_0 J_{3,27}-3 A_1 J_{6,27}-9 A_1 C_0
   J_{6,27}-3 A_0 C_1 J_{6,27}+3 A_2 J_{12,27}\\
   +9
   A_2 C_0 J_{12,27}+9 A_1 C_1
   J_{12,27}\big)
   +q^3 D \big(-3 A_1 J_{3,27}\\-9
   A_1 C_0 J_{3,27}-3 A_0 C_1 J_{3,27}-3 A_2
   J_{6,27}-9 A_2 C_0 J_{6,27}-9 A_1 C_1
   J_{6,27}+9 A_2 C_1 J_{12,27}\big)\\
   +q^4 D
   \left(-3 A_2 J_{3,27}-9 A_2 C_0 J_{3,27}-9 A_1
   C_1 J_{3,27}-9 A_2 C_1 J_{6,27}\right)-9 q^5 D
   A_2 C_1 J_{3,27}\\
   \equiv
 D\left[A_0 J_{12,27} +3 A_0 C_0  J_{12,27}-3 q^3 \left( A_0 C_1  J_{3,27}+A_1   J_{3,27}+ A_2  J_{6,27}\right)\right]\\
 +q D \left[-3 A_0 C_0    J_{6,27}+3 A_0 C_1  J_{12,27}-A_0    J_{6,27}+3 A_1  J_{12,27} -3  q^3  A_2 J_{3,27}\right]\\
  +q^2 D \left(-3 A_0 C_0  J_{3,27}-3 A_0 C_1  J_{6,27}-A_0  J_{3,27}-3 A_1  J_{6,27}+3 A_2    J_{12,27}\right)
 \pmod 9.
\end{multline}
By comparing powers of $q$ with exponent $\equiv 0 \pmod 3$, it may be seen that \eqref{conjmodn1eq3a} will hold if
\[
 3C_0 J_{12,27}- 3C_1 q^3 J_{3,27}\equiv 0 \pmod 9,
\]
or equivalently, if
\begin{align}\label{conjmod910mod3}
\sum_{n=0}^{\infty}(b_{3n}-a_{3n})q^{3n}
&\equiv DA_0(3C_0 J_{12,27}- 3C_1 q^3 J_{3,27})\pmod 9\\
&\equiv DA_0 \bigg( \left[2 \left(\frac{f_1^3}{f_3}-1\right) -3q\frac{f_9^3}{f_3}\right]J_{12,27}- 3 q^3 \frac{f_9^3}{f_3}J_{3,27}\bigg)\pmod 9\notag \\
&\equiv 0 \pmod 9, \notag
\end{align}
where the next-to-last congruence uses \eqref{3C03C1eq2}, and the last congruence follows from \eqref{JJ1}.

Similarly, by considering powers of $q$ with exponent $\equiv 1 \pmod{3}$, it follows that
\begin{align}\label{conjmod911mod3}
\sum_{n=0}^{\infty}(b_{3n+1}&+2a_{3n+1})q^{3n+1}\\
&\equiv -3 D q \left(A_0 C_0 J_{6,27}-A_0 C_1 J_{12,27}+3
   A_2 q^3 J_{3,27}+A_0 J_{6,27}-3 A_1
   J_{12,27}\right)\pmod 9\notag\\
   &\equiv -3DA_0q\left[\left( C_0+1\right) J_{6,27}- C_1 J_{12,27}\right]  \pmod 9\notag\\
   &\equiv -DA_0q\bigg(
\left[2 \left(\frac{f_1^3}{f_3}-1\right)-3q\frac{f_9^3}{f_3} +3\right]J_{6,27}- 3  \frac{f_9^3}{f_3}J_{12,27}\bigg)\notag\\
&\equiv 0 \pmod 9. \notag
\end{align}
The last two congruences follow from \eqref{3C03C1eq2} and \eqref{JJ2}, respectively.

Finally, by considering powers of $q$ with exponent $\equiv 2 \pmod 3$, it follows that
\begin{align}\label{conjmod912mod3}
\sum_{n=0}^{\infty}&(2b_{3n+2}+a_{3n+2})q^{3n+2}\\
&\equiv -3 D q^2 \left(2 A_0 C_0 J_{3,27}+2 A_0 C_1
   J_{6,27}+A_0 J_{3,27}+3 A_1 J_{6,27}-3 A_2
   J_{12,27}\right) \pmod 9 \notag\\
  & \equiv 3DA_0q^2\left[\left(C_0-1\right) J_{3,27}+C_1 J_{6,27}\right]\pmod 9 \notag\\
&\equiv Dq^2 \bigg(\left[2 \left(\frac{f_1^3}{f_3}-1\right)-3q\frac{f_9^3}{f_3} -3\right]J_{3,27}+ 3  \frac{f_9^3}{f_3}J_{6,27}  \bigg) \notag\\
&\equiv 0 \pmod 9, \notag
\end{align}
where once again \eqref{3C03C1eq2} has been used to derive the next-to-last congruence, and the final congruence follows from \eqref{JJ3}.
\end{proof}

\begin{corollary}\label{coranbn0mod9}
  Let the eta quotients $A(q)$ and $B(q)$ and the sequences $\{a_n\}$ and $\{b_n\}$ be as in Theorem \ref{conjmod91}. Then
  \begin{equation}\label{coranbn0mod9eq1}
   a_n \equiv 0 \pmod 9 \lrla b_n \equiv 0 \pmod 9.
  \end{equation}
\end{corollary}
\begin{proof}
  This is immediate from \eqref{conjmodn1eq3a} - \eqref{conjmodn1eq3c}.
\end{proof}

\begin{remark}
    One can see that $A(q)$ and $B(q)$ have identically vanishing coefficents modulo~9 non-trivially, since $f_{3}/f_{1}^{3}\not\equiv1\pmod{9}$.
\end{remark}

We next give some arithmetic/combinatorial consequences of Corollary \ref{coranbn0mod9}.
Recall that a  bipartition $\pi$ of the positive integer $n$ is a
pair of partitions $(\pi_1, \pi_2)$ with $|\pi_1| + |\pi_2| = n$, where as usual $|\lambda|$ denotes the sum of the parts of the partition $\lambda$. Let $p_{\_2}(n)$ denote the number of bipartitions $(\pi_1, \pi_2)$ of $n$
\begin{corollary}\label{corevodbipart}
  Let $S=\mathbb{N}\setminus 3\mathbb{N}$, the set of positive integers that are not multiples of 3, and let $D_s(n)$ denote the number of partitions of $n$ into an even number of distinct parts from $S$ minus the number of partitions of $n$ into an odd number of distinct parts from $S$. Then
  \begin{equation}\label{corevodbiparteq1}
    D_s(n)\equiv 0 \pmod 9 \lrla p_{\_2}(n) \equiv 0 \pmod 9.
  \end{equation}
\end{corollary}

\begin{proof}
   In Corollary \ref{coranbn0mod9}, set
 \begin{equation}\label{corevodbiparteq2}
  A(q) = \frac{f_1}{f_3}=(q,q^2;q^3)_{\infty}=\sum_{n=0}^{\infty}D_s(n) q^n\hspace{5pt} \lra \hspace{5pt} B(q)=A(q)\frac{f_3}{f_1^3}=\frac{1}{f_1^2}=\sum_{n=0}^{\infty}p_{\_2}(n)q^n.
  \end{equation}
  \end{proof}

Recall also that a partition is called \emph{3-regular} if none of its parts are multiples of 3.
\begin{corollary}
  Let $p_{\_2}^{(3)}(n)$ denote the number of bipartitions $(\pi_1, \pi_2)$ of $n$ where $\pi_1$ is 3-regular. Then
  \begin{equation}
  p_{\_2}^{(3)}(n)\equiv 0 \pmod 9 \lrla n \text{ is \underline{not} a generalized pentagonal number}.
\end{equation}
\end{corollary}
\begin{proof}
  In Corollary \ref{coranbn0mod9} let $A(q)=f_1$, so that
  \[
  B(q) = A(q) \frac{f_3}{f_1^3}= \frac{f_3}{f_1}\frac{1}{f_1}= \sum_{n=0}^{\infty}p_{\_2}^{(3)}(n)q^n,
  \]
  and recall the series expansion for $f_1$ at \eqref{f1ser}.
\end{proof}
Next, recall that when we say some property $P$ holds for ``almost all $n$'' in regard to a sequence $\{a_n\}$, we mean that
\[
\lim_{X\to \infty} \frac{\#\{n\leq X| P \text{ holds for }a_n \}}{X}=1.
\]
\begin{corollary}\label{corf17f3mod9}
 Let the sequences $\{c_n\}$ and $\{d_n\}$ be defined by
 \begin{equation}\label{corf17f3mod9eq1}
  f_1^7f_3 = \sum_{n=0}^{\infty}c_nq^n,  \hspace{25pt} \frac{f_1^7}{f_3} = \sum_{n=0}^{\infty}d_nq^n.
 \end{equation}
 Then $c_n \equiv 0 \pmod 9$ for almost all $n$, and $d_n \equiv 0 \pmod 9$ for almost all $n$.
\end{corollary}
\begin{proof}
For the first claim, let $A(q)=f_1^{10}$ in Corollary \ref{coranbn0mod9}, so that $B(q)=f_1^7f_3$ and recall that $f_1^{10}$ is lacunary. For the second claim, let $A(q) = f_1^7/f_3$, so that $B(q) = f_1^4$, recalling that $f_1^4$ is likewise lacunary.
\end{proof}

\begin{remark}
Note that by Serre's criteria for the vanishing of the coefficients in the series expansions of $f_1^{10}$  and $f_1^4$, we have that $c_n \equiv 0 \pmod 9$ if $12n+5$ has a prime factor $p\equiv-1\pmod{4}$ with odd exponent, and $d_n \equiv 0 \pmod 9$ if $6n+1$ has a prime factor $p\equiv-1\pmod{3}$ with odd exponent.
\end{remark}

\section{Vanishing Coefficients Modulo 25}

\subsection{Trivial Congruences.} Before considering congruence results that are not entirely trivial, we gain recall the two situations where the congruence results are trivial.

The first is where
\[
A(q)=: \sum_{n=0}^{\infty}a_nq^n, \hspace{25pt} B(q)=A(q)\frac{f_1^{25}}{f_5^5}=:  \sum_{n=0}^{\infty}b_nq^n,
\]
where $A(q)$ is any eta quotient. Then
\begin{equation}\label{zeromod25}
 a_n \equiv b_n \pmod{25}, \,\,\forall n \geq 0.
\end{equation}

The second trivial situation, not so easy to recognize, is where $A(q)$ and $B(q)$ have similar $m$-dissections for some positive integer $m$. We illustrate this with an example.

Let $C(q^4)=\prod_{j}f_j^{n_j}$ be any eta quotient such that $4|j$ for all $j$. Let
\[
A(q) = f_1^2f_2^3 C(q^4) =: \sum_{n=0}^{\infty}a_nq^n, \hspace{25pt} B(q)= \frac{f_1^6f_4^2}{f_2^3} C(q^4)=:  \sum_{n=0}^{\infty}b_nq^n.
\]
Then
\begin{equation}\label{zeroidmod25}
a_n \equiv 0 \pmod {25} \lrla b_n \equiv 0 \pmod {25}.
\end{equation}

This is trivially true once one knows that
\begin{align}
    f_1^2f_2^3&=  \frac{f_8^{15}}{f_4^4 f_{16}^6}-2 q\frac{ f_8^9}{f_4^2 f_{16}^2}-4 q^2 f_8^3
   f_{16}^2+8 q^3 \frac{f_4^2 f_{16}^6}{f_8^3}, \label{l5diss1b}\\
  \frac{f_1^6f_4^2}{f_2^3}&= \frac{f_8^{15}}{f_4^{4}f_{16}^6} -6 q \frac{f_8^{9}}{f_4^{2}f_{16}^2}
  +12 q^2 f_8^{3}f_{16}^2 - 8 q^3 \frac{f_{16}^6}{f_8^{3}}, \label{l5diss1c}
  \end{align}
  so that the $a_n$ and $b_n$ are non-zero multiples of each other, with the multipliers being relatively prime to 25.
  See  \cite[Lemma 2.7]{HMcLY23a}, for proofs of \eqref{l5diss1b} and \eqref{l5diss1c}.

  \subsection{The computer search for pairs of eta quotients with coefficients vanishing identically modulo 25.}

  To summarize, the situation that is first investigated in this section is the existence of pairs of eta quotients $(A(q), B(q))$ where
  \[
A(q)=: \sum_{n=0}^{\infty}a_nq^n, \hspace{25pt} B(q)=:  \sum_{n=0}^{\infty}b_nq^n, \hspace{25pt} 1 \leq j \leq 4,
\]
where
\begin{align}\label{zeromod252}
 &a_n \not \equiv b_n \pmod{25}, \,\,\forall n \geq 0,\\
 &a_n \equiv 0 \pmod {25} \lrla b_n \equiv 0 \pmod {25}. \notag
\end{align}

 We did not find any results that appeared to hold for infinite families of eta quotients, such as in Theorem \ref{conjmod91}.

  As regards the search for pairs of eta quotients ($A(q)$, $B(q)$) where \eqref{zeroidmod25} holds but \eqref{zeromod25} does not hold, we take our cue from Section \ref{mod9} and Equation \eqref{conjmod91eq2}, and mostly restrict our search to quintuples of eta quotients of the form
  \begin{equation}\label{FGfi5f5eq}
  \{F(q), G_1(q), G_2(q), G_3(q), G_4(q)\},\hspace{25pt} \text{ where }
  G_j(q) := \left(\frac{f_1^5}{f_5}\right)^j F(q), \hspace{25pt} 1 \leq j \leq 4.
  \end{equation}

Table \ref{ta5} (go to \href{https://www.wcupa.edu/sciences-mathematics/mathematics/jMcLaughlin/documents/identicalvanishingmod49and25tables.pdf}{https://tinyurl.com/529p5bjv}  to access the full version of the two tables, as the present paper contains just abbreviated versions of these) summarizes the results of some computer investigations. Columns 3 -7 show the counts of coefficients that are $0 \pmod{25}$ in the series expansion of $F(q)(f_1^5/f_5)^j$, $0\leq j \leq 4$, and columns 8 - 12 counts the number of zero coefficients in these same series expansions (high counts here indicating possible lacunarity of eta quotients, or coefficients possibly vanishing in arithmetic progressions).

As can be seen from this table, the most common situation may be described as follows. When one considers the set of five eta quotients described at \eqref{FGfi5f5eq} ($F(q)$ being the eta quotient listed in column 2 of Table \ref{ta5}), then four of the five eta quotients belong to a set all of whose coefficients vanish identically modulo 25, i.e., if $A(q)$ and $B(q)$ are any two eta quotients in this set of four, then \eqref{zeromod252} holds. If $C(q)=:\sum c_n q^n$ represents the fifth eta quotient in the set of five, and  $A(q)=:\sum a_n q^n$ is any of the other four, then there is strict inclusion of the sets of coefficients that vanish modulo 25, in that  $a_n \equiv 0 \pmod{25} \lra c_n \equiv 0 \pmod{25}$, and there exist integers $m$ such that $c_m \equiv 0 \pmod{25}$, but $a_m \not \equiv 0 \pmod{25}$.

We have included a number of rows to indicate the various kinds of behaviour that may occur. High counts of at least several thousand in one of columns 8 - 12 indicate that the corresponding eta quotient is likely to be lacunary. The most common situation is that one  of the five eta quotients represented in any particular row appears to be lacunary, while the other four are not. However, this is not universal, as the numbers in rows 34 and 48 indicate that three of the five eta quotients are lacunary. The numbers 1250, 1252 and 1253 in rows 49 and 50 may indicate coefficients that vanish in various arithmetic progressions, while the low numbers in rows 30 and 31 indicate that the eta quotients are likely to be neither lacunary nor to have coefficients that vanish in arithmetic progressions. The numbers 256 and 476 in rows 41 - 43 are possibly indicative of coefficients that vanish in arithmetic progressions, and we may investigate this phenomenon further in a subsequent paper.

Rows 252 - 254 are different from the other 251 rows in the table, as experiment suggests that in each case all five eta quotients have coefficients that vanish identically modulo 25.

 We give some proofs later to indicate how the various equalities and inclusions of sets of coefficients that vanish modulo 25, as suggested by the tables, may actually be demonstrated.


In Table \ref{ta6} below, the situation is somewhat different from that exhibited in Table \ref{ta5}. As in Table \ref{ta5}, Columns 3 -7 show the counts of coefficients that are $0 \pmod{25}$ in the series expansion of $F(q)(f_1^5/f_5)^j$, $0\leq j \leq 4$, and columns 8 - 12 counts the number of zero coefficients in these same series expansions (high counts here again indicating possible lacunarity of eta quotients, or coefficients possibly vanishing in arithmetic progressions).

This time, when one considers the set of five eta quotients described at \eqref{FGfi5f5eq} ($F(q)$ being the eta quotient listed in column 2 of Table \ref{ta6}), then just three of the five eta quotients belong to a set all of whose coefficients vanish identically modulo 25, i.e., if $A(q)$ and $B(q)$ are any two eta quotients in this set of three, then \eqref{zeromod252} holds.

One can see from columns 3 - 7 that there is a fourth eta quotient that has exactly one more vanishing coefficient modulo 25 in the first 15000 coefficients. Presumably there would be additional vanishing coefficients modulo 25 for this fourth eta quotient, compared to the other three, if we examined the coefficients further out than the first 15000.

If $C(q)=:\sum c_n q^n$ represents the fifth eta quotient in the set of five, and  $A(q)=:\sum a_n q^n$ is any of the other four, then there is strict inclusion of the sets of coefficients that vanish modulo 25, in that  $a_n \equiv 0 \pmod{25} \lra c_n \equiv 0 \pmod{25}$, and there exist integers $m$ such that $c_m \equiv 0 \pmod{25}$, but $a_m \not \equiv 0 \pmod{25}$.

We will illustrate the meaning of the numbers in columns 13 and 14 by considering row 1 in Table \ref{ta6} as an example.
{\allowdisplaybreaks
\begin{center}
\begin{longtable}{|C|C|CCCCC|CCCCC|CC|}
\hline
n&F(q)&0&1&2&3&4&0&1&2&3&4&C_N&N\\
\hline
1 & f_1 f_2 & 10505 & 7436 & 7436 & 7436 & 7437 & 10500 & 58 & 0 & 0 & 0 & 10441 & 5076 \\
\hline
\end{longtable}
\end{center}
}

It can be seen that the eta quotients represented by the  columns numbered 1, 2 and 3 (actually columns 4 - 6 in the table) each has 7436 coefficients that vanish modulo 25 (and in each case, vanish in the same 7436 place. The number 7437 in the column labelled 4 indicates that this eta quotient has one additional vanishing coefficient, modulo 25. The number 5076 in column 14 indicates that this additional coefficient is the 5076th out of 7437, and the 10441 in column 13 indicates that it is the coefficient of $q^{10441}$ in its series expansion.

It can again be seen that the most common situation appears to be  that at least one of the five eta quotients is lacunary, although row 8 in Table \ref{ta6} would seem to indicate that this is not universally so. As with Table \ref{ta5}, intermediate counts such as the numbers 58 and 474 in rows 1 and 3 respectively, may indicate coefficients that vanish in arithmetic progressions.

We indicate one curious phenomenon. The number 7141 that occurs in the $C_N$ column for rows 3, 6, 8, 9, 23, 28, 30, 44 and rows 46 - 49 indicates that the one additional zero coefficient mentioned above is in each case the coefficient of $q^{7141}$. There are additional cases of this phenomenon for exponents other than 7141 in the full version of the table. At present we do not have an explanation of this phenomenon.


\subsection{The case of  \texorpdfstring{$f_1^5f_5$}{} and \texorpdfstring{$f_1^{10}$}{}}
We return to row 48 in Table \ref{ta5}:
{\allowdisplaybreaks
\begin{center}
\begin{longtable}{|C|C|CCCCC|CCCCC|}
\hline
n&F(q)&0&1&2&3&4&0&1&2&3&4\\
\hline
48 & f_5^2 & 13693 & 7571 & 7571 & 7571 & 7571 & 13684 & 6123 & 6123 & 0 & 0 \\
\hline
\end{longtable}
\end{center}
}
With the notation of Equation \eqref{FGfi5f5eq}, we see that $G_1(q)=f_1^5f_5$ and $G_2(q)=f_1^{10}$, and that columns 4 and 5 above suggest that \eqref{zeromod252} holds, with $\{A(q),B(q)\} =\{f_1^{10},f_1^5f_5\}$.

We shall prove this by determining exactly when the coefficients in the series expansions for each of the two products vanish  modulo 25, and showing that in each case that the conditions are identical. As will be seen, this method of proof using modular forms is lengthy, so one might hope that a faster method might be found.

We also emphasize that the computer algebra system \emph{Mathematica} played a pretty much essential part in some of the proofs in this section. As an example, the relations stated in Lemma \ref{S3S4Hjl}
were first determined using \emph{Mathematica} (and subsequently verified by checking up to the Sturm bound). As a second example, the statements at \eqref{f110eq45} and \eqref{f110eq46} were generated using \emph{Mathematica}, and the eighteen other similar cases referred to but not displayed were similarly checked using \emph{Mathematica}. The reader will find many other instances in the various proofs where \emph{Mathematica} played a significant role.

We first consider $f_1^{10}$. Several intermediate lemmas are required.
After applying the dilation $q\to q^{12}$ and multiplying by $q^5$ to turn $f_1^{10}$ into a modular form, one gets
\begin{multline}\label{f110eq1}
 q^5 f_{12}^{10}  =q^5-10 q^{17}+35 q^{29}-30 q^{41}-105 q^{53}+238 q^{65}-260 q^{89}\\
 -165 q^{101}+140 q^{113}+1054
   q^{125}-770 q^{137}-595 q^{149}-715 q^{173}+2162 q^{185}+455 q^{197}+\dots .
\end{multline}
\begin{lemma}
One  has that
\begin{equation}\label{f110eq2}
  q^5 f_{12}^{10} = -\frac{1}{96}S_1(q)+\frac{1}{96}S_2(q),
\end{equation}
where $S_1(q)$ is the CM form of weight 5 and level 144 labelled \href{https://www.lmfdb.org/ModularForm/GL2/Q/holomorphic/144/5/g/a/}{144.5.g.a} in the \href{https://www.lmfdb.org/}{LMFDB}, and
$S_2(q)$ is the form labelled \href{https://www.lmfdb.org/ModularForm/GL2/Q/holomorphic/144/5/g/b/}{144.5.g.b}.
\end{lemma}
\begin{proof}
See Serre \cite[Eq. (27)]{S85}, or compare coefficients up to the Sturm bound.
\end{proof}
Since we need to give theta series representations for $S_1(q)$ and $S_2(q)$, for clarity of exposition, we state their initial series expansion:
\begin{multline}\label{f110eq3}
  S_1(q)= q-48 q^5+238 q^{13}+480 q^{17}+1679 q^{25}-1680 q^{29}+2162 q^{37}\\
  +1440 q^{41}+2401 q^{49}
  +5040   q^{53}   -6958 q^{61}-11424 q^{65}-1442 q^{73}-23040 q^{85}+12480 q^{89}\\
  +1918 q^{97}+7920
   q^{101}
   -9362 q^{109}   -6720 q^{113}+14641 q^{121}-50592 q^{125}\\+36960 q^{137}+80640 q^{145}+28560
   q^{149}
   -20398 q^{157}+28083 q^{169}+34320 q^{173}\\
   +64078 q^{181}-103776 q^{185}-38398
   q^{193}-21840 q^{197}+\dots,
\end{multline}
\begin{multline}\label{f110eq4}
   S_2(q) = q+48 q^5+238 q^{13}-480 q^{17}+1679 q^{25}+1680 q^{29}+2162 q^{37}\\
   -1440 q^{41}+2401 q^{49}-5040
   q^{53}-6958 q^{61}+11424 q^{65}-1442 q^{73}-23040 q^{85}-12480 q^{89}\\
   +1918 q^{97}-7920
   q^{101}-9362 q^{109}+6720 q^{113}+14641 q^{121}+50592 q^{125}\\
   -36960 q^{137}+80640 q^{145}-28560
   q^{149}-20398 q^{157}+28083 q^{169}-34320 q^{173}\\
   +64078 q^{181}+103776 q^{185}-38398
   q^{193}+21840 q^{197}+\dots .
\end{multline}

\begin{lemma}\label{apcongl}
The following identities hold.
  \begin{align}\label{SiHjeqs}
       S_{1}(q)&=H_{3}-H_{4}+iH_{7}-iH_{8},\\
    S_{2}(q)&=H_{3}-H_{4}-iH_{7}+iH_{8},
\end{align}
where
\begin{align}\label{Hjeqs}
     H_{3}&=\sum_{m,n}(6m+1+6ni)^{4}q^{(6m+1)^2+(6n)^2},\\
    H_{4}&=\sum_{m,n}(6m+3+(6n-2)i)^{4}q^{(6m+3)^2+(6n-2)^2},\\
    H_{7}&=\sum_{m,n}(6m+1+(6n-2)i)^{4}q^{(6m+1)^2+(6n-2)^2},\\
    H_{8}&=\sum_{m,n}(6m+1+(6n+2)i)^{4}q^{(6m+1)^2+(6n+2)^2}.
\end{align}
\end{lemma}
\begin{proof}
See \cite[Lemma 2.3, part (6)]{HMcLY22}, or again compare coefficients up to the Sturm bounds.
\end{proof}
Define the sequences $\{a_n\}$ and $\{b_n\}$ by
\[
S_1(q) =: \sum_{n=0}^{\infty}a_n q^n, \hspace{50pt} S_2(q) =: \sum_{n=0}^{\infty}b_n q^n.
\]
For later use, we summarize some of the properties of these coefficients.

\begin{lemma}\label{apcong2l}
The following hold:
\begin{itemize}
  \item If $n \not \equiv 1, 5 \pmod{12}$, $a_n=b_n =0$;
  \item If $n \equiv 1 \pmod{12}$, $a_n=b_n$;
  \item If $n \equiv 5 \pmod{12}$, $a_n=-b_n$;
  \item If $p \equiv 1 \pmod{12}$ is prime, with $p=x^2+y^2$ for integers $x>y>0$, then
  \begin{equation}\label{f110eq44a}
   a_p = \pm 2 \left(x^2-2 x y-y^2\right) \left(x^2+2 x y-y^2\right);
  \end{equation}
  \item If $p \equiv 5 \pmod{12}$ is prime, with $p=x^2+y^2$ for integers $x>y>0$, then
  \begin{equation}\label{f110eq44b}
   a_p = \pm 8 x y (x-y)(x+y);
  \end{equation}
  \item Apart from $p=5$, if $p \equiv 5 \pmod{12}$ is prime, $a_p \equiv 0 \pmod 5$;
   \item If $p \equiv 1 \pmod{12}$ is prime, $a_p  \equiv 2 \text{ or } 3 \pmod 5$;
  \item The recurrence formula at prime powers is
  \begin{equation}\label{f1f10eq5}
    a_{p^{k+1}}= a_{p^{k}}a_p-\chi(p) p^4  a_{p^{k-1}},
  \end{equation}
  where
  \begin{equation}\label{chipeq}
  \chi(p) =
  \begin{cases}
    1, & \mbox{if } p\equiv 1 \pmod 4 \\
     -1, & \mbox{if } p\equiv 3 \pmod 4.
  \end{cases}
  \end{equation}
\end{itemize}
\end{lemma}
\begin{proof}
These mostly follow from \eqref{SiHjeqs} and \eqref{Hjeqs}, after noting that the form of the general term in each of the $H_i$ is
\begin{equation}\label{f110eq6}
(x+i y)^4 q^{x^2+y^2}= \left[\left(x^2+2 x y-y^2\right) \left(x^2-2 x y-y^2\right)+4 i x y (x-y) (x+y) \right]q^{x^2+y^2}.
\end{equation}
The congruence statements for $a_p$, where $p>5$ is prime, $p=x^2+y^2$, make use of the facts that at most one of $x$, $y$ is $\equiv 0 \pmod 5$, and if neither is $\equiv 0 \pmod 5$, then
\[
\{x,y\}\pmod 5 \in \{\{1,1\}, \{1,4\},\{4,4\},\{2,2\},\{2,3\},\{3,3\} \}.
\]
The formula for $\chi(p)$ at \eqref{chipeq} may be accessed at the LMFDB page for the form labelled \href{https://www.lmfdb.org/ModularForm/GL2/Q/holomorphic/144/5/g/a/}{144.5.g.a}.
\end{proof}

Note for later use that the sequence $\{a_{p^n}\}$ is entirely determined by $a_p$, $\chi(p)$ and $p^4$, and hence are determined entirely modulo 5 and modulo 25 by these values.

We next consider when $a_{p^k}\equiv 0 \pmod{5}$.

\begin{lemma}\label{apk05l}
Let $p$ be an odd prime and $k$ a non-negative integer.\\
(1) If $p \equiv 3 \pmod 4$, $p \not =3$, then
\begin{equation}\label{p3mod4con}
a_{p^{2k+1}}=0, \hspace{25pt} a_{p^{2k}}=p^{4k}.
\end{equation}
(2) If $p=5$, then $a_{5^k} \equiv (a_5)^k \not \equiv 0 \pmod 5$.\\
(3) If $p \equiv 5 \pmod {12}$,   then
\[
a_{p^{2k+1}}\equiv 0 \pmod 5, \hspace{25pt} |a_{p^{2k}}|\equiv p^{4k} \not \equiv 0 \pmod 5.
\]
(4) If $p \equiv 1 \pmod {12}$,   then
\begin{equation}\label{f110eq9}
  a_{p^{5k+4}}\equiv 0 \pmod 5, k=0,1,2,\dots,
\end{equation}
and $a_{p^n}\not \equiv 0 \pmod 5$, if $n \not = 5k+4$, some non-negative integer $k$.
\end{lemma}
\begin{proof}
The proofs of (1), (2) and (3) follow directly from \eqref{f1f10eq5}, making use for (3) of the fact that if $p \equiv 5 \pmod {12}$, then $a_p \equiv 0 \pmod 5$.

If $p\equiv 1 \pmod{12}$ then from Lemma \ref{apcong2l} we have that $a_p\equiv 2,3 \pmod 5$. We consider each case in turn, making use of the fact that if $p\equiv 1 \pmod{12}$ is prime, then $p^4\equiv 1 \pmod 5$. If $a_p \equiv 2 \pmod 5$, then one obtains recursively from \eqref{f1f10eq5} that
\begin{equation}\label{f110eq7}
\{a_{p^0}, a_{p^1}, a_{p^2}, a_{p^3}, a_{p^4}, \dots \} \equiv \{1,2,3,4,0,1,2,3,4,0,\dots\}\pmod 5,
\end{equation}
with the indicated pattern repeating in steps of 5.   Likewise, if $a_p \equiv 3 \pmod 5$, then one obtains recursively that
\begin{equation}\label{f110eq8}
\{a_{p^0}, a_{p^1}, a_{p^2}, a_{p^3}, a_{p^4}, \dots \} \equiv \{1,3,3,1,0,4,2,2,4,0,1,3,3,1,0,4,2,2,4,0,\dots\}\pmod 5,
\end{equation}
with the indicated pattern repeating in steps of 10. This completes the proof for (4).
\end{proof}

To determine when $a_{p^n}\equiv 0 \pmod{25}$, we consider three cases
\begin{itemize}
  \item $a_p \equiv 0 \pmod{25}$,
  \item $a_p \equiv 0 \pmod 5, \text{ but } a_p \not \equiv 0 \pmod{25}$,
  \item  $a_p \not \equiv 0 \pmod 5$.
\end{itemize}

In what follows, the $x$ and $y$ will be the positive integers defined by $p=x^2+y^2$, for the prime $p$ (if relevant, we take $x>y$).

\begin{lemma}\label{apk025l}
Let $p\equiv 1 \pmod 4$ be a prime and $k$ a positive integer.\\
(1) If $a_p\equiv 0 \pmod{25}$, then
\begin{equation}\label{ap025eq1}
 a_{p^{2k+1}}\equiv 0 \pmod{25}, \hspace{25pt} |a_{p^{2k}}|\equiv p^{4k} \not \equiv 0 \pmod 5,
\end{equation}
with $a_p\equiv 0 \pmod{25}$ holding if 25 divides  exactly one of $x$, $y$, $x-y$ or $x+y$.\\
(2) If  $a_p \equiv 0 \pmod 5$ but $a_p \not \equiv 0 \pmod{25}$, then
\begin{equation}\label{f110eq11}
 a_{p^n}\equiv
 \begin{cases}
   0 \pmod{25}, & \mbox{if } n \equiv 9 \pmod{10}, \\
    0 \pmod{5}, \not \equiv 0 \pmod{25} & \mbox{if } n \equiv 1 \pmod{2},n \not \equiv 9 \pmod{10}, \\
     1, 2, 3 \text{ or } 4 \pmod{5}, & \mbox{if } n \equiv 0 \pmod{2}.
 \end{cases}
\end{equation}
(3) If $p \equiv 1 \pmod{12}$  then
\begin{equation}\label{f110eq17}
 a_{p^n}\equiv
 \begin{cases}
   0 \pmod{25}, & \mbox{if } n \equiv 24 \pmod{25}, \\
    0 \pmod{5}, \not \equiv 0 \pmod{25} & \mbox{if } n \equiv 4 \pmod{5}, n \not \equiv 24 \pmod{25}, \\
     1, 2, 3 \text{ or } 4 \pmod{5}, & \mbox{if } n \equiv 0,1,2,3 \pmod{5}.
 \end{cases}
\end{equation}
\end{lemma}

\begin{proof}
 The congruence statements at (1) follow from \eqref{f1f10eq5}, while the divisibility by 25 statement follows from \eqref{f110eq44b}, together with the fact that $\gcd(x,y)=1$.

 Next, if $a_p \equiv 0 \pmod 5$,  but $a_p \not \equiv 0 \pmod{25}$, there are four possibilities, namely, $a_p\equiv 5, 10, 15, 20 \pmod{25}$. Similarly, since $p^4\equiv 1 \pmod 5$, there are 5 possibilities, namely $p^4 \equiv 1, 6, 11, 16, 21 \pmod {25}$. Thus there are 20 possibilities to be considered. If, for example, $a_p \equiv 10 \pmod{25}$ and $p^4 \equiv 6 \pmod{25}$ then \eqref{f1f10eq5} gives that
\begin{multline}\label{f110eq10}
\{a_{p^0}, a_{p^1}, a_{p^2}, a_{p^3}, a_{p^4}, \dots \} \\
\equiv \{1,10,19,5,11,5,9,10,21,0,24,15,6,20,14,20,16,15,4,0,\\
\phantom{asdf}1,10,19,5,11,5,9,10,21,0,24,15,6,20,14,20,16,15,4,0,\\
1,10,19,5,11,5\dots\}\pmod{25},
\end{multline}
with the indicated pattern repeating in steps of 10, so that $a_{p^{10k+9}}\equiv 0 \pmod{25}$ for $k=0,1,2, \dots$, and otherwise the sequence also follows the pattern indicated at (2).

A similar situation holds for the other 19 cases (which can easily be checked with the assistance of a computer algebra system such as \emph{Mathematica}), so that if $p \equiv 5 \pmod{12}$ and $a_p \equiv 0 \pmod 5$ but $a_p \not \equiv 0 \pmod{25}$ then (2) holds.

Finally, we consider the case $a_p \not \equiv 0 \pmod 5$ (when $p \equiv 1 \pmod{12}$). As was seen above, in this case $a_p  \equiv 2, 3 \pmod 5$ or $a_p  \equiv 2, 3, 7, 8, 12, 13, 17, 18, 22, 23 \pmod{25}$. When this is combined with $p^4 \equiv 1, 6, 11, 16, 21 \pmod {25}$, it may seem that there should be 50 cases to consider. However, there are just 10 distinct cases modulo 25:
\begin{equation}\label{f110eq12}
(a_p,p^4)\pmod{25} \in \{ (2, 1), (3, 21), (7, 6), (8, 16), (12, 11), (13, 11), (17, 16), (18, 6), (22, 21), (23, 1)\}.
\end{equation}
One explanation for this is as follows. Recall that if $p=x^2+y^2$, then $a_p=\pm 2 \left(x^2-2 x y-y^2\right)$  $\left(x^2+2 x y-y^2\right)$. One can check, preferably using a computer algebra system, that for each
\begin{multline}\label{f110eq13}
  (x,y)\in \{(5, 1), (9, 1), (10, 1), (15, 1), (16, 1), (20, 1), (8, 2), (17, 2), (22, 3), (21, 4), (7, 5),\\ (18, 5), (24, 5), (11, 6), (14, 6), (10, 7), (12, 7), (13, 7), (15, 7), (20, 7), (23, 8), (24, 9), (18, 10), \\
   (24, 10), (19, 11), (18, 12), (18, 13), (19, 14), (18, 15), (24, 15), (24, 16), (23, 17), (20, 18), (24, 20)\}
\end{multline}
one gets that
\begin{equation}\label{f110eq14}
 (a_p,p^4)\pmod{25} = (2, 1).
\end{equation}

An exhaustive check, again using a computer algebra system, over all pairs $(x,y)$ with $0\leq y \leq y \leq 24$ with $x^2+y^2\not \equiv 0 \pmod 5$, leads to \eqref{f110eq12}.

Proceeding as above and using \eqref{f1f10eq5}, preferably in conjunction with a computer algebra system, that if, for example, \eqref{f110eq14} holds, then
\begin{multline}\label{f110eq15}
\{a_{p^0}, a_{p^1}, a_{p^2}, a_{p^3}, a_{p^4}, \dots \} \\
\equiv \{1,2,3,4,5,6,7,8,9,10,11,12,13,14,15,16,17,18,19,20,21,22,23,24,0,\\
1,2,3,4,5,6,7,8,9,10,11,12,13,14,15,16,17,18,19,20,21,22,23,24,0,\\
1,2,3,4,5,6,7,8,9,10,11\dots\}\pmod{25},
\end{multline}
where the indicated pattern repeats in steps of 25, showing that the statement at (3) holds. A similar situation holds for the other 9 cases in \eqref{f110eq12}. For example, if $(a_p,p^4)\pmod{25} = (3, 21)$, then
\begin{multline}\label{f110eq16}
\{a_{p^0}, a_{p^1}, a_{p^2}, a_{p^3}, a_{p^4}, \dots \} \\
\equiv \{1,3,13,1,5,19,2,7,4,15,11,18,23,16,15,9,12,22,14,5,21,8,8,6,0,\\
\phantom{as}24,22,12,24,20,6,23,18,21,10,14,7,2,9,10,16,13,3,11,20,4,17,17,19,0,\\
1,3,13,1,5,19,2,7,4,15,11\dots\}\pmod{25},
\end{multline}
where this time the pattern repeats in steps of 50, again showing that the statement at (3) holds, as it does for the other 8 cases.
\end{proof}

Upon all this together we get the following theorem on the vanishing modulo 25 of the coefficients in the series expansion of $f_1^{10}$.
\begin{theorem}
  Let the sequence $\{c_n\}$ be defined by
  \begin{equation}\label{f110eq18}
  f_1^{10}=:\sum_{n=0}^{\infty}c_nq^n.
  \end{equation}
  Then $c_n \equiv 0 \pmod{25}$ if and only if one of the following conditions hold:
  \begin{itemize}\label{f110eq19}
    \item $ord_p(12n+5)$ is odd for some prime $p \equiv -1 \pmod 4$;
    \item $ord_p(12n+5)$ is odd for some prime $p \equiv 5 \pmod {12}$ with $p=x^2+y^2$ such that 25 divides exactly one of $x$ $y$,  $x - y$ or $x + y$;
    \item $ord_p(12n+5) = 10k+9$, some integer $k\geq 0$, for some prime $p \equiv 5\pmod{12}$;
     \item $ord_p(12n+5) =25k+24$, some integer $k\geq 0$, for some prime $p \equiv 1\pmod{12}$;
     \item $12n+5$ is divisible by distinct primes $p_1, p_2 \equiv 1 \pmod 4$, such that
      \begin{itemize}
        \item if $p_1 \equiv 1 \pmod{12}$, then $ord_{p_1}(12n+5) =5k+4$, some integer $k\geq 0$ and $k \not \equiv 4 \pmod 5$;
        \item if $p_1 \equiv 5 \pmod{12}$, then $ord_{p_1}(12n+5)$ is odd, but not of the form $10k+9$;
        \item similar conditions hold for the prime $p_2$.
      \end{itemize}
  \end{itemize}
\end{theorem}
\begin{proof}
  These statements all follow from the fact that
  \[
  c_n = -\frac{1}{48}a_{12n+5},
  \]
 so that $c_n \equiv 0 \pmod{25} \lrla a_{12n+5}\equiv 0 \pmod{25}$, and the properties that have been proven for the coefficients $a_n$.
\end{proof}

We next show that if the sequence $\{d_n\}$ be defined by
  \begin{equation}\label{f110eq20}
  f_1^{5}f_5=:\sum_{n=0}^{\infty}d_nq^n,
  \end{equation}
 then $d_n \equiv 0 \pmod{25}$ under the exact same conditions. It will be seen that the proof is a little more technical, primarily due to the number of theta series in the linear combinations needed to express the CM forms $S_3(q)$, $\bar{S}_3(q)$, $S_4(q)$ and $\bar{S}_4(q)$ in  \eqref{f110eq21} below.

 After a similar dilation $q\to q^{12}$ and multiplying by $q^5$, again to produce the corresponding modular form, one gets that
 One has  the  series expansion
\begin{multline}\label{f110eq22}
  q^5 f_{12}^{5}f_{60}= q^5-5 q^{17}+5 q^{29}+10 q^{41}-15 q^{53}-7 q^{65}+20 q^{89}+5 q^{101}-5 q^{113}\\
  +14 q^{125}-35
   q^{137}-35 q^{149}+55 q^{173}+7 q^{185}+65 q^{197}+\dots
\end{multline}

 \begin{lemma}\label{f15f5l1}
  One has that
    \begin{equation}\label{f110eq21}
  q^5 f_{12}^{5}f_{60} = \left(-\frac{1}{32}+\frac{i}{24}\right)S_3(q)+\left(-\frac{1}{32}-\frac{i}{24}\right)\bar{S}_3(q)+
  \left(\frac{1}{32}-\frac{i}{24}\right)S_4(q)+  \left(\frac{1}{32}+  \frac{i}{24}\right)\bar{S}_4(q),
\end{equation}
where $S_3(q)$ is the CM form of weight 3 and level 720 labelled \href{https://www.lmfdb.org/ModularForm/GL2/Q/holomorphic/720/3/j/b/}{720.3.j.b} in the \href{https://www.lmfdb.org/}{LMFDB}, $\bar{S}_3(q)$ is its $i \to -i$ conjugate,
$S_4(q)$ is the form labelled \href{https://www.lmfdb.org/ModularForm/GL2/Q/holomorphic/720/3/j/c/}{720.3.j.c}, and $\bar{S}_4(q)$ is its $i \to -i$ conjugate.
 \end{lemma}
\begin{proof}
This follows upon comparing coefficients up to the Sturm bound.
\end{proof}

The form $S_3(q)$ has series expansion
\begin{multline}\label{f110eq23}
  S_3(q)= q-(4+3 i) q^5-24 i q^{13}+30 i q^{17}+(7+24 i) q^{25}-40 q^{29}+24 i q^{37}-80 q^{41}-49 q^{49}+90 i
   q^{53}+22 q^{61}\\
   -(72-96 i) q^{65}+96 i q^{73}+(90-120 i) q^{85}-160 q^{89}-144 i q^{97}-40
   q^{101}-182 q^{109}+30 i q^{113}+121 q^{121}\\
   +(44-117 i) q^{125}+210 i q^{137}+(160+120 i)
   q^{145}+280 q^{149}+264 i q^{157}-407 q^{169}-330 i q^{173}\\
   +38 q^{181}+(72-96 i) q^{185}-336 i
   q^{193}-390 i q^{197}+\dots,
\end{multline}
while
$S_4(q)$ has series expansion
\begin{multline}\label{f110eq24}
   S_4(q) = q+(4+3 i) q^5-24 i q^{13}-30 i q^{17}+(7+24 i) q^{25}+40 q^{29}+24 i q^{37}+80 q^{41}-49 q^{49}-90 i
   q^{53}+22 q^{61}\\
   +(72-96 i) q^{65}+96 i q^{73}+(90-120 i) q^{85}+160 q^{89}-144 i q^{97}+40
   q^{101}-182 q^{109}-30 i q^{113}+121 q^{121}\\
   -(44-117 i) q^{125}-210 i q^{137}+(160+120 i)
   q^{145}-280 q^{149}+264 i q^{157}-407 q^{169}+330 i q^{173}\\
   +38 q^{181}-(72-96 i) q^{185}-336 i
   q^{193}+390 i q^{197}+\dots .
\end{multline}

Before discussing the result in the next lemma, we introduce some notation to allow a collection of 96 theta series to be written efficiently. For $1 \leq j\leq 96$ let $(u_j(m), v_j(n))$ denote the $j$-th entry in the following list:
{\allowdisplaybreaks
\begin{multline}\label{uvlist}
\{(30 m-23,30 n+2),(30 m-23,30 n+4),(30 m-23,30 n+12),(30 m-23,30 n+14),\\
(30 m-23,30   n+22),(30 m-23,30 n+24),(30 m-21,30 n+2),(30 m-21,30 n+4),\\
   (30 m-21,30 n+14),(30 m-21,30  n+22),(30 m-13,30 n+2),(30 m-13,30 n+4),\\
   (30 m-13,30 n+12),(30 m-13,30 n+14),(30 m-13,30   n+22),(30 m-13,30 n+24),\\
   (30 m-11,30 n+2),(30 m-11,30 n+4),(30 m-11,30 n+12),(30 m-11,30  n+14),\\
   (30 m-11,30 n+22),(30 m-11,30 n+24),(30 m-3,30 n+2),(30 m-3,30 n+4),\\
   (30 m-3,30   n+14),(30 m-3,30 n+22),(30 m-1,30 n+2),(30 m-1,30 n+4),\\
   (30 m-1,30 n+12),(30 m-1,30   n+14),(30 m-1,30 n+22),(30 m-1,30 n+24),\\
   (30 m+1,30 n),(30 m+1,30 n+2),(30 m+1,30   n+4),(30 m+1,30 n+10),\\
   (30 m+1,30 n+12),(30 m+1,30 n+14),(30 m+1,30 n+20),(30 m+1,30   n+22),\\
   (30 m+1,30 n+24),(30 m+3,30 n+2),(30 m+3,30 n+4),(30 m+3,30 n+10),\\
   (30 m+3,30   n+14),(30 m+3,30 n+20),(30 m+3,30 n+22),(30 m+5,30 n+2),\\
   (30 m+5,30 n+4),(30 m+5,30   n+12),(30 m+5,30 n+14),(30 m+5,30 n+22),\\
   (30 m+5,30 n+24),(30 m+11,30 n),(30 m+11,30   n+2),(30 m+11,30 n+4),\\
   (30 m+11,30 n+10),(30 m+11,30 n+12),(30 m+11,30 n+14),(30 m+11,30   n+20),\\
   (30 m+11,30 n+22),(30 m+11,30 n+24),(30 m+13,30 n),(30 m+13,30 n+2),\\
   (30 m+13,30   n+4),(30 m+13,30 n+10),(30 m+13,30 n+12),(30 m+13,30 n+14),\\
   (30 m+13,30 n+20),(30   m+13,30 n+22),(30 m+13,30 n+24),(30 m+15,30 n+2),\\
   (30 m+15,30 n+4),(30 m+15,30 n+14),(30  m+15,30 n+22),(30 m+21,30 n+2),\\
   (30 m+21,30 n+4),(30 m+21,30 n+10),(30 m+21,30 n+14),(30   m+21,30 n+20),\\
   (30 m+21,30 n+22),(30 m+23,30 n),(30 m+23,30 n+2),(30 m+23,30 n+4),\\
   (30   m+23,30 n+10),(30 m+23,30 n+12),(30 m+23,30 n+14),(30 m+23,30 n+20),\\
   (30 m+23,30   n+22),(30 m+23,30 n+24),(30 m+25,30 n+2),(30 m+25,30 n+4),\\
   (30 m+25,30 n+12),(30 m+25,30   n+14),(30 m+25,30 n+22),(30 m+25,30 n+24)
\}
\end{multline}
}
For $1\leq j \leq 96$, let the Hecke theta series $H_j=H_j(q)$ be defined by
\begin{equation}\label{Hj96}
  H_j = \sum_{m,n=-\infty}^{\infty}(u_j(m)+i\,v_j(n))^2 q^{u_j(m)^2+v_j(n)^2}
\end{equation}
For $1 \leq j \leq 96$ let $\alpha_j$, respectively $\beta_j$, $\gamma_j$, $\delta_j$, denote the $j$ entry in, respectively, the following lists:
\begin{multline}\label{alphalist}
 \mathcal{A} =
 \{
 i,i,-1,-i,-i,1,0,-1,-1,0,-i,-i,-1,i,i,1,0,i,0,-i,0,1,1,-1,-1,1,0,-i,0,\\
 i,0,1,1,i,-i,i,-1,i,-i,-i,-1,-1,0,1,0,1,-1,i,i,1,-i,-i,-1,1,-i,i,-i,-1,\\
 -i,i,i,-1,-1,-i,0,-i,1,0,i,i,0,-1,1,1,-1,1,1,-1,1,-1,\\
 1,-1,i,0,i,1,0,-i,-i,0,-i,-i,1,i,i,-1
 \}
\end{multline}

\begin{multline}\label{betalist}
 \mathcal{B} =
 \{
 -i,0,1,0,i,0,1,1,1,1,i,0,1,0,-i,0,i,-i,-1,i,-i,-1,-1,0,0,-1,-i,i,-1,\\
 -i,i,-1,1,0,i,i,0,-i,-i,0,1,1,-1,1,-1,1,1,i,i,1,-i,-i,-1,1,0,-i,-i,0\\
 ,i,i,0,1,-1,i,i,-i,-1,-i,i,-i,1,-1,1,1,-1,0,-1,-1,-1,\\
 -1,0,-1,-i,-i,i,-1,i,-i,i,1,-i,-i,1,i,i,-1
 \}
\end{multline}

\begin{multline}\label{gammalist}
 \mathcal{G} =
 \{
 -i,-i,-1,i,i,1,0,-1,-1,0,i,i,-1,-i,-i,1,0,-i,0,i,0,1,1,-1,-1,1,0,i,0,\\
 -i,0,1,1,-i,i,-i,-1,-i,i,i,-1,-1,0,1,0,1,-1,-i,-i,1,i,i,-1,1,i,-i,i,\\
 -1,i,-i,-i,-1,-1,i,0,i,1,0,-i,-i,0,-1,1,1,-1,1,1,-1,\\
 1,-1,1,-1,-i,0,-i,1,0,i,i,0,i,i,1,-i,-i,-1
 \}
\end{multline}

\begin{multline}\label{deltalist}
 \mathcal{D} =
 \{
 i,0,1,0,-i,0,1,1,1,1,-i,0,1,0,i,0,-i,i,-1,-i,i,-1,-1,0,0,-1,i,-i,-1,i,\\
 -i,-1,1,0,-i,-i,0,i,i,0,1,1,-1,1,-1,1,1,-i,-i,1,i,i,-1,1,0,i,i,0,-i,-i,\\
 0,1,-1,-i,-i,i,-1,i,-i,i,1,-1,1,1,-1,0,-1,-1,-1,-1,0,\\
 -1,i,i,-i,-1,-i,i,-i,1,i,i,1,-i,-i,-1
 \}
\end{multline}

\begin{lemma}\label{S3S4Hjl}
 The following identities hold:
 \begin{align}\label{S3S4Hjeq}
   S_3(q) & =\sum_{j=1}^{96}\alpha_j H_j, \\
    \bar{S}_3(q) & =\sum_{j=1}^{96}\beta_j H_j, \\
     S_4(q) & =\sum_{j=1}^{96}\gamma_j H_j, \\
    \bar{S}_4(q) & =\sum_{j=1}^{96}\delta_j H_j.
 \end{align}
\end{lemma}
\begin{proof}
  These following upon comparing coefficients up to the Sturm bound.
\end{proof}

Define the sequences $\{e_n\}$ and $\{f_n\}$ by
\begin{equation}\label{S3S4cofseq}
S_3(q) =: \sum_{n=0}^{\infty}e_n q^n, \hspace{50pt} S_4(q) =: \sum_{n=0}^{\infty}f_n q^n.
\end{equation}

\begin{lemma}\label{enfnrelsl}
 The following hold:
\begin{itemize}
  \item If $n \not \equiv 1, 5 \pmod{12}$, $e_n=f_n =0$;
  \item If $n \equiv 1 \pmod{12}$, $e_n=f_n$;
  \item If $n \equiv 5 \pmod{12}$, $e_n=-f_n$;
  \item If $p>5$, $p \equiv 1 \pmod{4}$ is prime, with $p=x^2+y^2$ for integers $x>y>0$, then
  \begin{equation}\label{f110eq33a}
   e_p =
   \begin{cases}
     \pm 2(x^2-y^2), & \mbox{if } p \equiv 1,49 \pmod{60}, \\
      \pm 4 i x y, & \mbox{if } p \equiv 13, 37 \pmod{60}, \\
      \pm 2 i (x^2-y^2), & \mbox{if } p \equiv 17,53 \pmod{60}, \\
     \pm 4 x y, & \mbox{if } p \equiv 29,41 \pmod{60}.
   \end{cases}
  \end{equation}
  \item Apart from $p=5$, if $p \equiv 5 \pmod{12}$ is prime, $e_p \equiv 0 \pmod 5$;
   \item If $p \equiv 1 \pmod{12}$ is prime, $e_p \not \equiv 0 \pmod 5$;
  \item The recurrence formula at prime powers is
  \begin{equation}\label{f1f10eq35}
    e_{p^{k+1}}= e_{p^{k}}e_p-\chi(p) p^2  e_{p^{k-1}},
  \end{equation}
  where
  \[
  \chi(p) =
  \begin{cases}
    1, & \mbox{if } p\equiv 1,3,7,9 \pmod{20}\\
     -1, & \mbox{if } p\equiv 11,13,17,19 \pmod{20} .
  \end{cases}
  \]
\end{itemize}
\end{lemma}

\begin{proof}
The proof involves quite a lot of tedious examination of cases, preferably checked using a computer algebra system
For example, it can be checked that, modulo 60,
the exponents in all 96 theta series, lie in the set
\[
\{1, 5, 13, 17, 25, 29, 37, 41, 49, 53\},
\]
leading to a proof of the first statement.

Similarly, one finds that the exponent of the theta series $H_j$, namely $u_j(m)^2+v_j(n)^2$, is $\equiv 1 \pmod{12}$ for $j$ in the set
\begin{multline*}
  \{
  3, 6, 7, 8, 9, 10, 13, 16, 19, 22, 23, 24, 25, 26, 29, 32, 33, 37,
41, 42, 43, 44, 45, 46, 47, 50, \\
53, 54, 58, 62, 63, 67, 71, 72, 73,
74, 75, 76, 77, 78, 79, 80, 81, 82, 86, 90, 93, 96
  \},
\end{multline*}
and that for these values of $j$ one has $\alpha_j = \gamma_j$, leading to a proof of the second assertion.

Likewise, one has that the exponent of the theta series $H_j$,  $u_j(m)^2+v_j(n)^2$, is $\equiv 5 \pmod{12}$ for $j$ in the set
\begin{multline*}
  \{
 1, 2, 4, 5, 11, 12, 14, 15, 17, 18, 20, 21, 27, 28, 30, 31, 34, 35,
36, 38, 39, 40, 48, 49, \\51, 52, 55, 56, 57, 59, 60, 61, 64, 65, 66,
68, 69, 70, 83, 84, 85, 87, 88, 89, 91, 92, 94, 95
  \},
\end{multline*}
and that for these values of $j$ one has $\alpha_j = -\gamma_j$, giving a proof of the third assertion.

For the fourth item, notice that the general term in each theta series has the form
\[
(x+ i\, y)^2 q^{x^2+y^2} = [(x^2-y^2)+2i\,x y]q^{x^2+y^2}.
\]
If $p\equiv 1 \pmod 4$ is a prime, then the representation $p=x^2+y^2$ in integers $x$ and $y$ is unique up to sign and interchanging $x$ and $y$. By considering the exponents of $q$ in all 96 $H_j$, taking into account the forms of the $(u_j(m), v_j(n))$ at \eqref{uvlist}, it can be shown that if $p$ is represented the exponent in some theta series, then it has exactly two representations, either both coming from the same theta series, or coming from two different theta series. If one term that contributes to the coefficient of $q^p = q^{x^2+y^2}$ is $[(x^2-y^2)+2i\,x y]q^{x^2+y^2}$ then the other term that contributes has the same form but with exactly one of $x$ or $y$ replaced with its negative, and thus has the form $[(x^2-y^2)-2i\,x y]q^{x^2+y^2}$.

If $p$ has two representations coming from the same theta series, then either $p \equiv 1 \pmod{60}$ or $p \equiv 49 \pmod{60}$, with the corresponding $(u_j(m),v_j(n))$ in the former case having one of three forms, namely,
\[
(30m+1, 30n), \hspace{25pt} (30m+15,30n+4),  \hspace{25pt}  (30m+15, 30n+14),
\]
while also having one of three forms in the later case, namely,
\[
(30m+13, 30n), \hspace{25pt} (30m+15,30n+2),  \hspace{25pt}  (30m+15, 30n+22).
\]
In all cases, the corresponding $\alpha_j$ are either 1 or -1, so that
\[
e_p = \pm \bigg([(x^2-y^2)+2i\,x y]+ [(x^2-y^2)-2i\,x y] \bigg)=\pm 2(x^2-y^2).
\]
If $p$ is represented by two different theta series, say $H_j$ and $H_k$, then one finds (once again assisted by a computer algebra system)
that there four possibilities for the corresponding pairs $(\alpha_j, \alpha_k)$:
\[
(\alpha_j, \alpha_k)=
\begin{cases}
 \pm (1,1), & \mbox{if } p \equiv 1, 49 \pmod{60}, \\
 \pm (1,-1), & \mbox{if } p \equiv 13, 37 \pmod{60}, \\
 \pm (i,i), & \mbox{if } p \equiv 17, 53 \pmod{60}, \\
 \pm (i,-i), & \mbox{if } p \equiv 29, 41 \pmod{60}. \\
\end{cases}
\]
The combinations
\[
\alpha_j[(x^2-y^2)+2i\,x y]+ \alpha_k[(x^2-y^2)-2i\,x y]
\]
then lead to the four statements in item 4.

The fifth assertion follows in the case that $p \equiv 17, 53 \pmod{60}$ from the facts that $e_p = \pm 2 i\,(x^2-y^2)=\pm 2(x-y)(x+y)$ and that from
\eqref{uvlist} the form of the pairs $\{x,y\}=\{u_j(m),v_j(n)\}$ that represent $p=x^2+y^2$, lie in the collection
\begin{multline*}
 \{\{30 m-11,30 n+4\},\{30 m-11,30 n+14\},\{30 m-1,30 n+4\},\{30 m-1,30 n+14\},\\
 \{30 m+1,30 n+4\},\{30
   m+1,30 n+14\},\{30 m+11,30 n+4\},\{30 m+11,30 n+14\}\}
\end{multline*}
for $p \equiv 17\pmod{60}$, and lie in the collection
\begin{multline*}
 \{\{30 m-23,30 n+2\},\{30 m-23,30 n+22\},\{30 m-13,30 n+2\},\{30 m-13,30 n+22\},\\
 \{30 m+13,30
   n+2\},\{30 m+13,30 n+22\},\{30 m+23,30 n+2\},\{30 m+23,30 n+22\}\}
\end{multline*}
for $p \equiv 53 \pmod{60}$.

When  $p \equiv 29, 41 \pmod{60}$ and $e_p = \pm 4xy$, the assertion follows form the facts that the pairs $\{x,y\}=\{u_j(m),v_j(n)\}$ that represent $p=x^2+y^2$, lie in the collection
\begin{multline*}
 \{\{30 m+5,30 n+2\},\{30 m+5,30 n+22\},\{30 m+13,30 n+10\},\{30 m+13,30
   n+20\},\\
   \{30 m+23,30 n+10\},\{30 m+23,30 n+20\},\{30 m+25,30 n+2\},\{30 m+25,30
   n+22\}\}
\end{multline*}
for $p \equiv 29\pmod{60}$, and lie in the collection
\begin{multline*}
 \{\{30 m+1,30 n+10\},\{30 m+1,30 n+20\},\{30 m+5,30 n+4\},\{30 m+5,30
   n+14\},\\
   \{30 m+11,30 n+10\},\{30 m+11,30 n+20\},\{30 m+25,30 n+4\},\{30 m+25,30
   n+14\}\}
\end{multline*}
for $p \equiv 41 \pmod{60}$.

The sixth assertion follows from similar arguments. With the preceding notation, when $p \equiv 1, 49 \pmod{60}$  and $e_p = \pm 2 \,(x^2-y^2)=\pm 2(x-y)(x+y)$, the form of the pairs $\{x,y\}=\{u_j(m),v_j(n)\}$ that represent $p=x^2+y^2$ lie in the collection
\begin{multline*}
 \{
 \{30 m+1,30 n\},\{30 m+5,30 n+24\},\{30 m+11,30 n\},\{30 m+15,30 n+4\},\\
 \{30
   m+15,30 n+14\},\{30 m+21,30 n+10\},\{30 m+21,30 n+20\},\{30 m+25,30 n+24\}
   \}
\end{multline*}
for $p \equiv 1\pmod{60}$, and lie in the collection
\begin{multline*}
 \{
\{30 m+3,30 n+10\},\{30 m+3,30 n+20\},\{30 m+5,30 n+12\},\{30 m+13,30 n\},\\
\{30
   m+15,30 n+2\},\{30 m+15,30 n+22\},\{30 m+23,30 n\},\{30 m+25,30 n+12\}
   \}
\end{multline*}
for $p \equiv 49 \pmod{60}$, and it is easily seen that in all cases $5 \nmid (x^2-y^2)$. Likewise, when $p \equiv 13, 37 \pmod{60}$  and $e_p = \pm 4\,i\,xy$, the form of the pairs $\{x,y\}=\{u_j(m),v_j(n)\}$ that represent $p=x^2+y^2$ lie in the collection
\begin{multline*}
 \{
\{30 m-23,30 n+12\},\{30 m-13,30 n+12\},\{30 m-3,30 n+2\},\{30 m-3,30   n+22\},\\
   \{30 m+3,30 n+2\},\{30 m+3,30 n+22\},\{30 m+13,30 n+12\},\{30 m+23,30   n+12\}
   \}
\end{multline*}
for $p \equiv 13\pmod{60}$, and lie in the collection
\begin{multline*}
 \{
\{30 m-21,30 n+4\},\{30 m-21,30 n+14\},\{30 m-11,30 n+24\},\{30 m-1,30   n+24\},\\
   \{30 m+1,30 n+24\},\{30 m+11,30 n+24\},\{30 m+21,30 n+4\},\{30 m+21,30   n+14\}
   \}
\end{multline*}
for $p \equiv 37 \pmod{60}$, and it is easily seen that in all cases $5 \nmid 4 \, i \, xy$.

Finally, the formula for $\chi(p)$ at \eqref{f1f10eq35} may be found at the LMFDB page for newform \href{https://www.lmfdb.org/ModularForm/GL2/Q/holomorphic/720/3/j/b/}{720.3.j.b}.
\end{proof}

As we did when examining $f_1^{10}$, we next consider when $e_{p^k}\equiv 0 \pmod{5}$.

\begin{lemma}
(1)  If $p \equiv 3 \pmod 4$, $p \not =3$, then
\begin{equation}\label{f15f5p3congeq}
e_{p^{2k+1}}=0, \hspace{25pt} |e_{p^{2k}}|=p^{2k}.
\end{equation}
(2) If $p=5$, then $e_{5^k} \equiv (e_5)^k \not \equiv 0 \pmod 5$.\\
(3) If $p \equiv 5 \pmod {12}$,  then
\begin{equation}\label{f15f55congeq}
 e_{p^{2k+1}}\equiv 0 \pmod 5, \hspace{25pt} |e_{p^{2k}}|\equiv p^{2k} \not \equiv 0 \pmod 5.
\end{equation}
(4) If $p\equiv 1 \pmod{12}$,
\begin{equation}\label{f110eq39}
  e_{p^{5k+4}}\equiv 0 \pmod 5, k=0,1,2,\dots,
\end{equation}
and $e_{p^n}\not \equiv 0 \pmod 5$, if $n \not = 5k+4$, some non-negative integer $k$.
\end{lemma}
\begin{proof}
The statements at (1), (2) and (3) follow from \eqref{f1f10eq35}, here also using in the case of (3) that if $p \equiv 5 \pmod {12}$, then $e_p \equiv 0 \pmod 5$.

Next, consider $p\equiv 1 \pmod{12}$, prime, and $p=x^2+y^2$. If $p\equiv 1 \pmod 5$ ($p\equiv 1 \pmod{60}$) then
\[
(x,y)\pmod 5 \in \{(0, 1), (0, 4), (1, 0), (4, 0) \}\lra (p^2,e_p)\pmod 5 \in\{(1,2),(1,3)\}.
\]
Likewise, if  $p\equiv 4 \pmod 5$ ($p\equiv 49 \pmod{60}$) then
\[
(x,y)\pmod 5 \in \{(0, 2), (0, 3), (2, 0), (3, 0) \}\lra (p^2,e_p)\pmod 5 \in\{(1,2),(1,3)\}.
\]
If $p\equiv 3 \pmod 5$ ($p\equiv 13 \pmod{60}$) then
\[
(x,y)\pmod 5 \in \{(2, 2), (2, 3), (3, 2), (3, 3) \}\lra (p^2,e_p)\pmod 5 \in\{(4,i),(4,-i)\}.
\]
Finally, if  $p\equiv 2 \pmod 5$ ($p\equiv 37 \pmod{60}$) then
\[
(x,y)\pmod 5 \in \{(1, 1), (1, 4), (4, 1), (4, 4) \}\lra (p^2,e_p)\pmod 5 \in\{(4,i),(4,-i)\}.
\]
If one then uses \eqref{f1f10eq35} modulo 5 with
\begin{equation}\label{f110eq366}
(p^2,e_p)\pmod 5 \in\{(1,2),(1,3),(4,i),(4,-i)\}
\end{equation}
one gets that $\{e_{p^0}, e_{p^1}, e_{p^2}, e_{p^3}, e_{p^4}, \dots \}$ is congruent modulo 5, respectively, to
\begin{align}\label{f110eq37}
& \{1,2,3,4,0,1,2,3,4,0,1,2,3,4,0,1,2,3,4,0,\dots\},\\
& \{1,3,3,1,0,4,2,2,4,0,1,3,3,1,0,4,2,2,4,0,\dots\},\\
& \{1,i,3,2 i,0,-2 i,2,-i,4,0,1,i,3,2 i,0,-2 i,2,-i,4,0,\dots\},\\
& \{1,-i,3,-2 i,0,2 i,2,i,4,0,1,-i,3,-2 i,0,2 i,2,i,4,0,\dots\}.
\end{align}
with the indicated patterns repeating in steps of either 5 or 10.  Thus we have in all cases that (4) holds.

\end{proof}

As in the case of $f_1^{10}$, to determine when $e_{p^n}\equiv 0 \pmod{25}$, we consider three cases
\begin{itemize}
  \item $e_p \equiv 0 \pmod{25}$,
  \item $e_p \equiv 0 \pmod 5, \text{ but } e_p \not \equiv 0 \pmod{25}$,
  \item  $e_p \not \equiv 0 \pmod 5$.
\end{itemize}

\begin{lemma}\label{f15f5cong25l}
(1) If $e_p \equiv 0 \pmod{25}$, then
\begin{equation}\label{ep025eq1}
 e_{p^{2k+1}}\equiv 0 \pmod{25}, \hspace{25pt} |e_{p^{2k}}|\equiv p^{2k} \not \equiv 0 \pmod 5,
\end{equation}
and $e_p \equiv 0 \pmod{25}$ holds only when 25 divides exactly one of $x$, $y$, $x-y$ or $x+y$.\\
(2) If $p \equiv 5 \pmod{12}$ and $e_p \equiv 0 \pmod 5$ but $e_p \not \equiv 0 \pmod{25}$ then
\begin{equation}\label{f110eq41}
 e_{p^n}\equiv
 \begin{cases}
   0 \pmod{25}, & \mbox{if } n \equiv 9 \pmod{10}, \\
    0 \pmod{5}, \not \equiv 0 \pmod{25} & \mbox{if } n \equiv 1 \pmod{2},n \not \equiv 9 \pmod{10}, \\
     1, 2, 3 \text{ or } 4 \pmod{5}, & \mbox{if } n \equiv 0 \pmod{2}.
 \end{cases}
\end{equation}
(3) If $p \equiv 1 \pmod{12}$  then
\begin{equation}\label{f110eq47}
e_{p^n}\equiv
 \begin{cases}
   0 \pmod{25}, & \mbox{if } n \equiv 24 \pmod{25}, \\
    0 \pmod{5}, \not \equiv 0 \pmod{25} & \mbox{if } n \equiv 4 \pmod{5}, n \not \equiv 24 \pmod{25}, \\
     \text{not } 0 \pmod{5}, & \mbox{if } n \equiv 0,1,2,3 \pmod{5}.
 \end{cases}
\end{equation}
\end{lemma}

\begin{proof}
The congruence \eqref{ep025eq1} follows from \eqref{f1f10eq35}, where as earlier the latter non-congruence follows from the consideration of $e_p \equiv 0 \pmod{5}$ above. Note that if $e_p\equiv 0 \pmod{25}$, then $p \equiv 5 \pmod 12$, and thus either
\begin{itemize}
  \item $p\equiv 17 \text{ or } 53 \pmod{60} \lra e_p=\pm 2(x^2-y^2)\lra 25|(x-y) \text{ or } 25|(x+y)$,
  \item $p\equiv 29 \text{ or } 41 \pmod{60} \lra e_p=\pm 4xy\lra 25|x \text{ or } 25|y$.
\end{itemize}
Once again we have used the fact that $\gcd(x,y)=1$.
In addition, an exhaustive search modulo 60 shows that if $25|x$ (or $25|y$) and $p=x^2+y^2$ is prime, then $p \equiv 1, 29, 41 \text{ or } 49 \pmod{60}$, so that if in addition $p\equiv 5 \pmod{12}$, then $p \equiv  29 \text{ or }  41  \pmod{60}$. Likewise, if $25|(x-y)$ (or $25|(x-y)$) a similar exhaustive search modulo 60 shows that  $p \equiv 13, 17, 37 \text{ or } 53 \pmod{60}$, so that if in addition $p\equiv 5 \pmod{12}$, then $p \equiv  17 \text{ or }  53  \pmod{60}$. Thus  we have shown that $e_p\equiv 0 \pmod{25}$ exactly when $p \equiv 5 \pmod{12}$ and  25 divides exactly one of $x$, $y$, $x-y$ or $x+y$, where $p=x^2+y^2$.

Next, we consider $e_p \equiv 0 \pmod 5$,  but $e_p \not \equiv 0 \pmod{25}$. From \eqref{f110eq33a} if $p\equiv 17, 53 \pmod{60}$, then $e_p \equiv -10i,-5i, 5i \text{ or } 10i \pmod{25}$ and (since $p \equiv 2 \text{ or } 3 \pmod 5$) $p^2 \equiv 4,9,14, 19 \text{ or } 24 \pmod{25}$, thus giving 20 possibilities for the pair $(e_p, p^2)\pmod{25}$. Similarly, if $p\equiv 29, 41 \pmod{60}$,  then $e_p\equiv 5, 10, 15, 20 \pmod{25}$, and (since $p \equiv 2 \text{ or } 3 \pmod 5$) $p^2 \equiv 1,6,11, 16 \text{ or } 21 \pmod{25}$,  giving a further 20 possibilities for the pair $(e_p, p^2)\pmod{25}$.  Thus in total there are 40 possibilities to be considered. If $e_p \equiv -5i \pmod{25}$ and $p^2 \equiv 14 \pmod{25}$ then \eqref{f1f10eq35} gives that
\begin{multline}\label{f110eq40}
\{e_{p^0}, e_{p^1}, e_{p^2}, e_{p^3}, e_{p^4}, \dots \} \\
\equiv \{1,-5 i,14,10 i,21,10 i,19,-5 i,16,0,24,5 i,11,-10 i,4,-10 i,6,5 i,9,0,\\
\phantom{adf}1,-5 i,14,10 i,21,10 i,19,-5 i,16,0,24,5 i,11,-10 i,4,-10 i,6,5 i,9,0,\\
1,-5 i,14,10 i,21,10 i,\dots\}\pmod{25},
\end{multline}
with the indicated pattern repeating in steps of 10, so that $e_{p^{10k+9}}\equiv 0 \pmod{25}$ for $k=0,1,2, \dots$. As a second example, if If $e_p \equiv 15 \pmod{25}$ and $p^2 \equiv 11 \pmod{25}$ then \eqref{f1f10eq35} likewise gives that
\begin{multline}\label{f110eq40a}
\{e_{p^0}, e_{p^1}, e_{p^2}, e_{p^3}, e_{p^4}, \dots \} \\
\equiv \{1,15,14,20,21,20,19,15,16,0,24,10,11,5,4,5,6,10,9,0,\\
\phantom{adf}1,15,14,20,21,20,19,15,16,0,24,10,11,5,4,5,6,10,9,0,\\
1,15,14,20,21,20,\dots\}\pmod{25},
\end{multline}
with the indicated pattern also repeating in steps of 10, so that once again $e_{p^{10k+9}}\equiv 0 \pmod{25}$ for $k=0,1,2, \dots$.
A similar situation holds for the other 38 cases, so that \eqref{f110eq41} holds.

Finally, we consider the case $e_p \not \equiv 0 \pmod 5$ (when $p \equiv 1 \pmod{12}$). From \eqref{f110eq366}, either \\
   $e_p  \equiv 2, 3, 7, 8, 12, 13, 17, 18, 22, 23 \pmod{25}$ and $p^2 \equiv 1, 6, 11, 16, 21 \pmod {25}$, or\\
 $e_p  \equiv -i, i, -4 i, 4 i, -6 i, 6 i, -9 i, 9 i, -11 i, 11 i \pmod{25}$ and $p^2 \equiv 4, 9, 14, 19, 24 \pmod {25}$. It may seem that there should be 100 distinct cases to consider. However, as was the case when considering $f_1^{10}$, there  are a good many less, and indeed there are just 20 distinct cases modulo 25:
\begin{multline}\label{f110eq42}
(e_p,p^2)\pmod{25}\\
 \in \{
(-i, 19), (i, 19), (-4 i, 4), (4 i, 4), (-6 i, 9), (6 i, 9), (-9 i, 14), (9 i, 14), (-11 i, 24), (11 i, 24),\\
 (2, 1), (3, 21), (7, 6), (8, 16), (12, 11), (13, 11), (17, 16), (18, 6), (22, 21), (23, 1)
\}.
\end{multline}
Here also \eqref{f110eq42} follows from an exhaustive check  using a computer algebra system, over all pairs $(x,y)$ with $0\leq y \leq y \leq 24$ with $x^2+y^2\not \equiv 0 \pmod 5$. The behaviour modulo 25 in the sequence $e_{p^0}, e_{p^1}, e_{p^2}, e_{p^3}, e_{p^4}, \dots$ is similar in all cases. We give two examples.

If $(e_p,p^2)\pmod{25} = (i, 19)$ then \eqref{f1f10eq35} gives(with the assistance of a computer algebra system) that
\begin{multline}\label{f110eq45}
\{e_{p^0}, e_{p^1}, e_{p^2}, e_{p^3}, e_{p^4}, \dots \}\equiv \\
 \{1,i,18,12 i,5,8 i,12,-11 i,14,5 i,11,6 i,3,-8 i,15,-12 i,22,-6 i,24,10 i,21,11 i,13,-3 i,0,\\
-7 i,7,-i,9,-10 i,6,-9 i,23,2 i,10,-2 i,17,4 i,19,-5 i,16,-4 i,8,7 i,20,3 i,2,9 i,4,0,\\
1,i,18,12 i,5,8 i,\dots\}\pmod{25},
\end{multline}
where the indicated pattern repeats in steps of 50. A similar situation holds for the other 19 cases in \eqref{f110eq42}. For example, if $(e_p,p^2)\pmod{25} = (17, 16)$, then
\begin{multline}\label{f110eq46}
\{e_{p^0}, e_{p^1}, e_{p^2}, e_{p^3}, e_{p^4}, \dots \} \\
\equiv \{1,17,23,19,5,6,22,3,24,10,11,2,8,4,15,16,7,13,9,20,21,12,18,14,0,\\
\phantom{asd}1,17,23,19,5,6,22,3,24,10,11,2,8,4,15,16,7,13,9,20,21,12,18,14,0,\\
1,17,23,19,5,6,\dots\}\pmod{25},
\end{multline}
where this time the pattern repeats in steps of 25. Thus we have that \eqref{f110eq47} holds.

\end{proof}

Collecting these results together we get the following theorem, a partner to Theorem \ref{f110eq18}.
\begin{theorem}
  Let the sequence $\{d_n\}$ be defined by
  \begin{equation}\label{f110eq48}
  f_1^{5}f_5=:\sum_{n=0}^{\infty}d_nq^n.
  \end{equation}
  Then $d_n \equiv 0 \pmod{25}$ if and only if one of the following conditions hold:
  \begin{itemize}\label{f110eq49}
    \item $ord_p(12n+5)$ is odd for some prime $p \equiv -1 \pmod 4$;
    \item $ord_p(12n+5)$ is odd for some prime $p \equiv 1 \pmod 4$ with $p=x^2+y^2$ such that $25|x  (x - y)  y  (x + y)$;
    \item $ord_p(12n+5) = 10k+9$, some integer $k\geq 0$, for some prime $p \equiv 5\pmod{12}$;
     \item $ord_p(12n+5) =25k+24$, some integer $k\geq 0$, for some prime $p \equiv 1\pmod{12}$;
     \item $12n+5$ is divisible by distinct primes $p_1, p_2 \equiv 1 \pmod 4$, such that
      \begin{itemize}
        \item if $p_1 \equiv 1 \pmod{12}$, then $ord_{p_1}(12n+5) =5k+4$, some integer $k\geq 0$ and $k \not \equiv 4 \pmod 5$;
        \item if $p_1 \equiv 5 \pmod{12}$, then $ord_{p_1}(12n+5)$ is odd, but not of the form $10k+9$;
        \item similar conditions hold for the prime $p_2$.
      \end{itemize}
  \end{itemize}
\end{theorem}
\begin{proof}
Upon using \eqref{f110eq21} together with the relations that have been shown between the $e_n$ and $f_n$ in Lemma \ref{enfnrelsl}, one gets that if $e_{12n+5} = u+i\,v$, then
\begin{equation}\label{dnformeq}
 d_n = -\frac{u}{8}-\frac{v}{6}.
\end{equation}

If $5 \nmid 12n+5$, then  $e_{12n+5}$ is either purely real ($=u$, say) or purely imaginary ($=i\,v$,say), so that in this case
\[
d_n =-\frac{e_{12n+5}}{8}\text{ or } d_n =\frac{i\,e_{12n+5}}{6} \lra d_n \equiv 0 \pmod{25} \lrla e_{12n+5} \equiv 0 \pmod{25}.
\]
If $5|12n+5$ so that $12n+5 = 5^k m$ for some positive integer $k$ and $\gcd(5,m)=1$, then by multiplicativity

\[
e_{12n+5}= e_{5^k}e_m \equiv (e_5)^k e_m \equiv (1+2\,i)^k \pmod 5,
\]
with $e_m$ being either purely real ($=u$, say) or purely imaginary ($=i\,v$, say). If $(1+2\,i)^k \pmod 5=x + iy$,  then from \eqref{dnformeq}
\[
d_n \equiv \left( -\frac{x}{8}-\frac{y}{6}\right)u \pmod 5 \text{ or } d_n \equiv \left( \frac{y}{8}-\frac{x}{6}\right)v \pmod 5.
\]
Since
\[
\{(1+2\,i)^m :m \in \mathbb{N}\}\pmod 5 = \{1+2 i,2-i,-1-2 i,-2+i \},
\]
and it is an easy check that if $x+i\,y$ is any of the four numbers in the second collection that neither the numerator or denominator of either of
\[
\left( -\frac{x}{8}-\frac{y}{6}\right) \text{ or }\left( \frac{y}{8}-\frac{x}{6}\right)
\]
is divisible by 5.
Thus, when $e_{12n+5}= e_{5^k}e_m$ for some integer $k\geq 1$, $\gcd(5,m)=1$,
\[
25|d_n \lrla 25 |e_m \lrla 25|e_{12n+5}.
\]
Thus in all cases we have that
\[
d_n \equiv 0 \pmod{25} \lrla e_{12n+5} \equiv 0 \pmod{25},
\]
and the rest of the proof follows from what has been proven for the sequence $e_n$.
\end{proof}

Upon comparing Theorem \ref{f110eq18} and Theorem \ref{f110eq48}, we get the following result.
\begin{theorem}\label{cd25}
   Let the sequences  $\{c_n\}$  and  $\{d_n\}$ be defined by
  \begin{equation}\label{f110eq50}
  f_1^{10}=:\sum_{n=0}^{\infty}c_nq^n, \hspace{25pt}  f_1^{5}f_5=:\sum_{n=0}^{\infty}d_nq^n.
  \end{equation}
  Then
  \begin{equation}\label{f110eq51}
    c_n \equiv 0 \pmod{25} \lrla d_n \equiv 0 \pmod{25}.
  \end{equation}
\end{theorem}

\vspace{10pt}

\subsection{The case of  \texorpdfstring{$f_{1}f_{5}$}{} and  \texorpdfstring{$f_{1}^{6}$}{}}

In this subsection we consider an example of strict inclusion of vanishing modulo~25 between two eta quotients suggested by the row 34 of Table~\ref{ta5}.
{\allowdisplaybreaks
\begin{center}
\begin{longtable}{|C|C|CCCCC|CCCCC|}
\hline
n&F(q)&0&1&2&3&4&0&1&2&3&4\\
\hline
34 & f_1 f_5 & 12168 & 9207 & 9207 & 9207 & 9207 & 12161 & 7887 & 5661 & 0 & 0 \\
\hline
\end{longtable}
\end{center}
}

Looking at that line, one can note that if one defines
\begin{equation}\label{34ab}
\sum_{n=0}^{\infty}a_{n}q^{n}=f_{1}f_{5}\quad\mbox{and}\quad \sum_{n=0}^{\infty}b_{n}q^{n}=f_{1}^{6},
\end{equation}
then one is supposed to have
\begin{theorem}\label{abrow34}
    Let $a_{n}$ and $b_{n}$ be defined as in~\eqref{34ab}. Then
    $$
    \{n|\,b_{n}\equiv0\pmod{25}\}\subsetneqq \{n|\,a_{n}\equiv0\pmod{25}\}.
    $$
\end{theorem}

The approach to validate Theorem~\ref{abrow34} is similar to that of Theorem~\ref{cd25} and by characterizing the vanishing of the eta quotients modulo~25 via the CM representations for $qf_{4}f_{20}$ and $qf_{4}^{6}$ given in \cite{HLMcLYYZ}. To that end, we first define
\begin{equation}
    \label{ABrow34}
    \sum_{n=0}^{\infty}A_{n}q^{n}=qf_{4}f_{20}\quad\mbox{and}\quad \sum_{n=0}^{\infty}B_{n}q^{n}=qf_{4}^{6}
\end{equation}
and note that $A_{4n+1}=a_{n}$, $B_{4n+1}=b_{n}$, and $A_{n}=B_{n}=0$ for $n\not\equiv1\pmod{4}$, so that Theorem~\ref{abrow34} amounts to
\begin{theorem}\label{ABrow34thm}
    Let $A_{n}$ and $B_{n}$ be defined as in~\eqref{ABrow34}. Then
    $$
    \{n|\,B_{n}\equiv0\pmod{25}\}\subsetneqq \{n|\,A_{n}\equiv0\pmod{25}\}.
    $$
\end{theorem}
By \cite{HLMcLYYZ}, it is known that both $A_{n}$ and $B_{n}$ are multiplicative, i.e.,
$$
A_{n}=\prod_{p|n}A_{p^{m}}\quad\mbox{and}\quad B_{n}=\prod_{p|n}B_{p^{m}}
$$
for $n=\prod_{p|n}p^{e_{p}}$ with $p^{e_{p}}||n$. So the vanishing of $A_{n}$ and $B_{n}$ modulo~25 can be described by that of the local factors $A_{p^{e_{p}}}$, $B_{p^{e_{p}}}$. As such, we first recall by \cite{HLMcLYYZ} a formula for $A_{p^{e_{p}}}$:

\begin{lemma}
    \label{prop1-34}
    Let $A_{n}$ be defined as in~\eqref{ABrow34}. Then one has that
    \begin{enumerate}
\item $A_{5^{m}}=(-1)^{m}$,

\item for $p\equiv 1,9\pmod{20}$,
 $$
A_{p^{m}}=\begin{cases}
   m+1 &\mbox{if $p=X^{2}+5Y^{2}$ with $2|Y$, }\\
  \left(-1\right)^{m}\left(m+1\right)  &\mbox{otherwise,}
\end{cases}
$$
    \item for $p\equiv3,7\pmod{20}$,
$$
A_{p^{m}}=\begin{cases}
    0&\mbox{if $m$ is odd,}\\
    \left(-1\right)^{m/2}&\mbox{otherwise,}
\end{cases}
$$
\item for $p\equiv 11,13,17,19\pmod{20}$,
$$
A_{p^{m}}=\begin{cases}
0&\mbox{if $m$ is odd,}\\
1&\mbox{otherwise.}
\end{cases}
$$
\end{enumerate}
\end{lemma}

An immediate implication of Lemma~\ref{prop1-34} is the following characterization of the vanishing of $A_{n}$ modulo~25.

\begin{proposition}
        \label{cor1-34}
    Let $n$ be a positive integer. Then $A_{n}\equiv0\pmod{25}$ if and only if one of the following holds true:
    \begin{enumerate}
        \item $n$ has a prime factor $p\equiv1,9\pmod{20}$ with exponent $e_{p}\equiv-1\pmod{25}$,
        \item $n$ has two distinct prime factors $p_{1},p_{2}\equiv1,9\pmod{20}$ with exponent $e_{p_{i}}\equiv-1\pmod{5}$,
        \item $n$ has a prime factor $p\equiv3,7,11,13,17,19\pmod{20}$ with odd exponent.
    \end{enumerate}
\end{proposition}

In what follows we shall make use of the CM representation for $qf_{4}^{6}$ as a holomorphic modular form to determine a necessary condition on a positive integer $n$ such that $B_{n}\equiv0\pmod{25}$. Such a CM representation can be found in \cite{HLMcLYYZ} and is stated in Lemma~\ref{prop2-34}. 

\begin{lemma}
    \label{prop2-34}
The following identity holds.
    $$
   qf_{4}^{6} =\frac{1}{2}\sum_{m,n=-\infty}^{\infty}\left(2m+1+2ni\right)^2q^{\mathcal{N}\left(2m+1+2ni\right)}.
$$
\end{lemma}

Using Lemma~\ref{prop2-34} one can formulate the residue of $B_{p^{m}}$ modulo~$5$ as follows.

\begin{lemma}\label{prop3-34}
Let $B_{n}$ be defined as in~\eqref{ABrow34}. Then the following hold true.
    \begin{enumerate}
        \item
        $$
  B_{5^{m}}\equiv \frac{1}{2}\left((1+2i)^{2m}+(1-2i)^{2m}\right)\ne0\pmod{5\mathbb{Z}[i]}.
  $$

  \item For $p\equiv 1,9\pmod{20}$,
$$
  B_{p^{m}}\equiv \frac{(\pm1)^{m}}{2}(m+1)\quad\mbox{or}\quad \frac{(\pm2)^{m}}{2}(m+1)\pmod{5\mathbb{Z}[i]}.
  $$

  \item For $p\equiv 3, 7, 11,19\pmod{20}$, $B_{p^{m}}=0$ for $m$ odd, and
  $$
  B_{p^{m}}=p^{m}\ne0\pmod{5}
  $$
  for $m$ even.

\item For $p\equiv 13,17\pmod{20}$,
  $$
  b(p^{m})\equiv \frac{1}{2}(1+(-1)^{m})(2i)^{m}\pmod{5\mathbb{Z}[i]},
  $$
  or
  $$
  b(p^{m})\equiv \frac{1}{2}(1+(-1)^{m})(4i)^{m}\pmod{5\mathbb{Z}[i]}.
  $$
    \end{enumerate}
\end{lemma}


\begin{proof}
  By the theta representation for $qf_{4}^{6}$ and basic knowledge of algebraic theory of quadratic orders, it is not hard to see that
  $$
  B_{5^{m}}=\frac{1}{2}\sum_{r=0}^{m}(1+2i)^{2r}(1-2i)^{2(m-r)}.
  $$
  Since $5=(1+2i)(1-2i)$, and $(1+2i)$ and $(1-2i)$ are coprime over $\mathbb{Z}[i]$, then moduloing $5\mathbb{Z}[i]$ one finds that
  $$
  B_{5^{m}}\equiv \frac{1}{2}\left((1+2i)^{2m}+(1-2i)^{2m}\right)\ne0\pmod{5\mathbb{Z}[i]}.
  $$

  For $p\equiv 1,9\pmod{20}$, it is clear that $p=(x+yi)(x-yi)$ for some $x$ odd and $y$ even such that if $p\equiv1\pmod{20}$ (resp. $p\equiv9\pmod{20}$) exactly one of $x$ and $y$ is congruent to~$\pm1$ (resp. $\pm2$) modulo~5, and the other is divisible by~5. Suppose that $x\equiv\pm1\pmod{5}$ and $5|y$. So one can deduce that
  $$
  B_{p^{m}}=\frac{1}{2}\sum_{r=0}^{m}(x+yi)^{2r}(x-yi)^{2(m-r)}\equiv \frac{(\pm1)^{m}}{2}(m+1)\pmod{5\mathbb{Z}[i]}.
  $$
  Similarly, if $x\equiv\pm2\pmod{20}$, then
   $$
  B_{p^{m}}\equiv \frac{(\pm2)^{m}}{2}(m+1)\pmod{5\mathbb{Z}[i]}.
  $$

  For $p\equiv 3, 7, 11,19\pmod{20}$,  since $p$ is inert in $\mathbb{Z}[i]$, then it is clear that $B_{p^{m}}=0$ for $m$ odd, and
  $$
  B_{p^{m}}=p^{m}\ne0\pmod{5}
  $$
  for $m$ even.

  For $p\equiv 13,17\pmod{20}$, so that $p\equiv2,3\pmod{5}$, it is easy to see that $p=x^{2}+y^{2}=(x+yi)(x-yi)$ with $x\equiv y\equiv\pm1\pmod{5}$ and $x\equiv y\equiv\pm2\pmod{5}$, respectively. So one can deduce that for $x\equiv y\equiv\pm1\pmod{5}$,
  $$
  B_{p^{m}}=\frac{1}{2}\sum_{r=0}^{m}(x+yi)^{2r}(x-yi)^{2(m-r)}\equiv \frac{1}{2}((\pm(1+i))^{2m}+(\pm(1-i))^{2m})\pmod{5\mathbb{Z}[i]},
  $$
  and thus,
  $$
  B_{p^{m}}\equiv \frac{1}{2}(1+(-1)^{m})(2i)^{m}\pmod{5\mathbb{Z}[i]}.
  $$
  Similarly, for $x\equiv y\equiv2\pmod{5}$, one has that
  $$
  B_{p^{m}}\equiv \frac{1}{2}(1+(-1)^{m})(4i)^{m}\pmod{5\mathbb{Z}[i]}.
  $$
\end{proof}

As an application of Lemma~\ref{prop3-34}, one can deduce necessary conditions on a positive integer $n$ for which $B_{n}\equiv0\pmod{25}$.

\begin{proposition}
    \label{cor2-34}
Let $B_{n}$ be defined as in~\eqref{ABrow34}. Then $B_{n}\equiv0\pmod{25}$ only if one of the following holds true:
    \begin{enumerate}
    \item $n$ has a prime factor $p\equiv1,9\pmod{20}$ with exponent $e_{p}\equiv-1\pmod{5}$,

    \item $n$ has a prime factor $p\equiv3,7,11,13,17,19\pmod{20}$ with odd exponent.

    \end{enumerate}
\end{proposition}

Finally, one has that

\begin{proof}[Proof of Theorem~\ref{ABrow34thm}]
    This follows from Propositions~\eqref{cor1-34} and~\eqref{cor2-34}.
\end{proof}

\section{Concluding Remarks}
There are several questions that can be asked, upon consideration of the results in this paper.

\begin{itemize}
  \item Are there similar results modulo $p^2$, for $p\geq 7$?
  \item Are there any general results for $p=5$ similar to those contained in Theorem \ref{mod4t1} and Theorem \ref{conjmod2nn}
for $p=2$ and Theorem \ref{conjmod91} for $p=3$?
\item Are there other combinatorial applications of any of the theorems, similar to those given in Examples \ref{ex1} and \ref{ex2} and Corollaries \ref{cex1} and \ref{oddpartscong}?
\item Are there combinatorial proofs of any of the combinatorial results stated?
\item In the case of $p=5$, are there methods of proof that do not involve modular forms?
\end{itemize}
One reason for the last question is that the methods we used apply only to \emph{lacunary} eta quotients, whereas as the tables would appear to indicate that in most cases where coefficients vanish identically modulo 25, the  eta quotients are \underline{not} lacunary. Also, as has been seen, the proofs we gave that used modular forms were quite technical, and it would clearly be advantageous to have simpler methods of proof.

Recall the function $D_S(n)$ from Corollary \ref{corevodbipart}, where
   $S$ denotes the set of positive integers that are not multiples of 3, and  $D_S(n)$  is the number of partitions of $n$ into an even number of distinct parts from $S$ minus the number of partitions of $n$ into an odd number of distinct parts from $S$.
   From \eqref{f13diss} one has that
   \begin{multline}\label{3dissid1}
 \sum_{n=0}^{\infty}D_s(n) q^n =(q,q^2;q^3)_{\infty}=    \frac{f_1}{f_3} = \frac{J_{12,27}}{f_3}-q \frac{J_{6,27}}{f_3}-q^2 \frac{J_{3,27}}{f_3}=\\
   \frac{1}{(q^3,q^6,q^9,q^{18},q^{21},q^{24};q^{27})_{\infty}}
-\frac{q}{(q^3,q^9,q^{12},q^{15},q^{18},q^{24};q^{27})_{\infty}}
-\frac{q^2}{(q^6,q^9,q^{12},q^{15},q^{18},q^{21};q^{27})_{\infty}}.
   \end{multline}
This dissection explains why
  $D_S(3n) \to \infty$ and  $D_S(3n+1), \, D_S(3n+12) \to -\infty$ as $n \to \infty$.
Note that this behaviour is very different from what happens if we replace $S$ with $\mathbb{N}$ and let $D_{\mathbb{N}}(n)$  be the number of partitions of $n$ into an even number of distinct parts from $\mathbb{N}$ minus the number of partitions of $n$ into an odd number of distinct parts from $\mathbb{N}$, when we recall Franklin's proof of the pentagonal number theorem \cite{F}, which showed
\[
\{D_{\mathbb{N}}(n)|n \in \mathbb{N}\} = \{-1,0,1\},
\]
with $D_{\mathbb{N}}(n)=0$ if $n$ is not a generalized pentagonal number. Moreover, \eqref{3dissid1} implies some interesting partition identities (after making the replacement $q\to q^{1/3}$ in the infinite products on the right side). For $a\in\{1,2,4\}$ let $p_{a,9}(n)$ denote the number of partitions of $n$ into parts $\not \equiv \pm a,0 \pmod 9$. Then
\begin{align}\label{3dissid2}
 D_S(3n) &= p_{4,9}(n), \\
  D_S(3n+1) &= -p_{2,9}(n),\notag \\
  D_S(3n+2) &= -p_{1,9}(n).\notag
\end{align}
As an example, if we take $n=20$, then the number of partitions of 60 into an \underline{even} number of distinct parts from $S$ is
631, the number of partitions of 60 into an \underline{odd} number of distinct parts from $S$ is 407, so that $D_{S}(60) = 631-407=224$, and one similarly computes that $p_{4,9}(20)=224$.

Are there combinatorial proofs of these identities?

\section{Appendix: Tables}

The full versions of the tables below may be found at \href{https://www.wcupa.edu/sciences-mathematics/mathematics/jMcLaughlin/documents/identicalvanishingmod49and25tables.pdf}{https://tinyurl.com/529p5bjv}, as we deemed them too long to include in the printed version of the paper.
  {\allowdisplaybreaks
\begin{center}
\LTcapwidth=0.8\textwidth
\begin{longtable}{|C|C|CCCCC|CCCCC|}
\caption{The count of zero coefficients modulo 25 (columns 3 -7) and the count of zero coefficients (columns 8 - 12) in the first 15000 terms in the series expansion of $F(q)(f_1^5/f_5)^j$, $0\leq j \leq 4$.}\label{ta5}\\
\hline
n&F(q)&0&1&2&3&4&0&1&2&3&4\\
\hline
 1 & f_1 & 14810 & 3199 & 3199 & 3199 & 3199 & 14800 & 22 & 0 & 0 & 0 \\
 \vdots &  \vdots &  \vdots&  \vdots &  \vdots &  \vdots &  \vdots &  \vdots &  \vdots &  \vdots &  \vdots & \vdots \\
 28 & \frac{f_3 f_4^2}{f_2} & 11968 & 3346 & 3346 & 3346 & 3346 & 11961 & 78 & 0 & 0 & 0 \\
 29 & \frac{1}{f_1 f_5} & 4715 & 4715 & 4715 & 5426 & 4715 & 0 & 0 & 2051 & 2 & 0 \\
 30 & \frac{f_1 f_2^9}{f_4 f_5} & 6581 & 8633 & 6581 & 6581 & 6581 & 3 & 1 & 0 & 0 & 0 \\
 31 & \frac{f_4^2}{f_1^2 f_2 f_5} & 2741 & 2741 & 2741 & 3454 & 2741 & 0 & 1 & 0 & 3 & 0 \\
 32 & f_5 & 14921 & 3068 & 3068 & 3068 & 3068 & 14911 & 23 & 0 & 0 & 0 \\
 33 & \frac{f_5}{f_1} & 6804 & 6804 & 6804 & 8570 & 6804 & 0 & 6408 & 0 & 4 & 0 \\
 34 & f_1 f_5 & 12168 & 9207 & 9207 & 9207 & 9207 & 12161 & 7887 & 5661 & 0 & 0 \\
  \vdots &  \vdots &  \vdots&  \vdots &  \vdots &  \vdots &  \vdots &  \vdots &  \vdots &  \vdots &  \vdots & \vdots \\
 41 & \frac{f_2^2 f_5}{f_4} & 11599 & 3865 & 3865 & 3865 & 3865 & 11591 & 256 & 0 & 0 & 0 \\
 42 & \frac{f_2^3 f_5}{f_1 f_4} & 11806 & 5071 & 5071 & 5071 & 5071 & 11797 & 476 & 0 & 0 & 0 \\
 43 & f_4 f_5 & 11411 & 5209 & 5209 & 5209 & 5209 & 11404 & 476 & 0 & 0 & 0 \\
 \vdots &  \vdots &  \vdots&  \vdots &  \vdots &  \vdots &  \vdots &  \vdots &  \vdots &  \vdots &  \vdots & \vdots \\
 48 & f_5^2 & 13693 & 7571 & 7571 & 7571 & 7571 & 13684 & 6123 & 6123 & 0 & 0 \\
 49 & \frac{f_5^4}{f_1} & 2490 & 2490 & 2490 & 3476 & 2490 & 0 & 1252 & 1 & 1250 & 0 \\
 50 & \frac{f_4^2 f_5^4}{f_1^2 f_2} & 2693 & 2693 & 2693 & 3617 & 2693 & 0 & 1250 & 2 & 1253 & 0 \\
  \vdots &  \vdots &  \vdots&  \vdots &  \vdots &  \vdots &  \vdots &  \vdots &  \vdots &  \vdots &  \vdots & \vdots \\
 81 & \frac{f_5 f_6}{f_1^2 f_2 f_3} & 1944 & 8371 & 1944 & 1944 & 1944 & 0 & 8365 & 2 & 0 & 0 \\
  \vdots &  \vdots &  \vdots&  \vdots &  \vdots &  \vdots &  \vdots &  \vdots &  \vdots &  \vdots &  \vdots & \vdots \\
 252 & \frac{f_4^2 f_5^2}{f_1 f_2} & 1257 & 1257 & 1257 & 1257 & 1257 & 0 & 135 & 0 & 0 & 0 \\
 253 & f_1 f_3^4 f_5 & 1188 & 1188 & 1188 & 1188 & 1188 & 477 & 0 & 1 & 0 & 0 \\
 254 & \frac{f_1^2 f_2^5 f_5}{f_4^2} & 1258 & 1258 & 1258 & 1258 & 1258 & 134 & 0 & 0 & 0 & 0 \\
 \hline
\end{longtable}
\end{center}
}

\footnotesize

  {\allowdisplaybreaks
\begin{center}
\LTcapwidth=0.8\textwidth
\begin{longtable}{C|C|CCCCC|CCCCC|CC}
\caption{The count of zero coefficients modulo 25 (columns 3 -7) and the count of zero coefficients (columns 8 - 12) in the first 15000 terms in the series expansion of $F(q)(f_1^5/f_5)^j$, $0\leq j \leq 4$. This time only 3 of the 5 eta quotients have identically vanishing coefficients, modulo 25.}\label{ta6}\\
n&F(q)&0&1&2&3&4&0&1&2&3&4&C_N&N\\
\hline
1 & f_1 f_2 & 10505 & 7436 & 7436 & 7436 & 7437 & 10500 & 58 & 0 & 0 & 0 & 10441 & 5076 \\
 \vdots &  \vdots  &  \vdots  &  \vdots  &  \vdots  &  \vdots  &  \vdots  &  \vdots  &  \vdots  &  \vdots  &  \vdots  &  \vdots  & \vdots  &  \vdots  \\
 3 & f_3^2 & 12881 & 3176 & 3177 & 3176 & 3176 & 12874 & 474 & 0 & 0 & 0 & 7141 & 1421 \\
 \vdots &  \vdots  &  \vdots  &  \vdots  &  \vdots  &  \vdots  &  \vdots  &  \vdots  &  \vdots  &  \vdots  &  \vdots  &  \vdots  & \vdots  &  \vdots  \\
 6 & \frac{f_1 f_3 f_4}{f_2} & 10639 & 2030 & 2031 & 2030 & 2030 & 10634 & 3 & 0 & 0 & 0 & 7141 &
   929 \\
 \vdots &  \vdots  &  \vdots  &  \vdots  &  \vdots  &  \vdots  &  \vdots  &  \vdots  &  \vdots  &  \vdots  &  \vdots  &  \vdots  & \vdots  &  \vdots  \\
 8 & \frac{f_1 f_2^{38}}{f_4^{14} f_5^3} & 1596 & 1596 & 7891 & 1597 & 1596 & 0 & 0 & 0 & 0 & 0 &
   7141 & 692 \\
 9 & \frac{f_5^3}{f_1 f_2^4} & 6120 & 6119 & 6119 & 7891 & 6119 & 0 & 0 & 5661 & 7887 & 0 & 7141 &
   2780 \\
 \vdots &  \vdots  &  \vdots  &  \vdots  &  \vdots  &  \vdots  &  \vdots  &  \vdots  &  \vdots  &  \vdots  &  \vdots  &  \vdots  & \vdots  &  \vdots  \\
 23 & \frac{f_4^2 f_6}{f_8} & 12812 & 2570 & 2570 & 2570 & 2571 & 12805 & 1 & 0 & 0 & 0 & 7141 &
   1215 \\
  \vdots &  \vdots  &  \vdots  &  \vdots  &  \vdots  &  \vdots  &  \vdots  &  \vdots  &  \vdots  &  \vdots  &  \vdots  &  \vdots  & \vdots  &  \vdots  \\
 28 & \frac{f_2 f_5^3}{f_1 f_{10}} & 7891 & 6119 & 6120 & 6119 & 6119 & 7887 & 5661 & 0 & 0 & 0 &
   7141 & 2780 \\
 \vdots &  \vdots  &  \vdots  &  \vdots  &  \vdots  &  \vdots  &  \vdots  &  \vdots  &  \vdots  &  \vdots  &  \vdots  &  \vdots  & \vdots  &  \vdots  \\
 30 & \frac{f_1^2 f_6^4}{f_2 f_3^2 f_{12}} & 10639 & 2150 & 2151 & 2150 & 2150 & 10634 & 1 & 0 & 0 &
   0 & 7141 & 1026 \\
 \vdots &  \vdots  &  \vdots  &  \vdots  &  \vdots  &  \vdots  &  \vdots  &  \vdots  &  \vdots  &  \vdots  &  \vdots  &  \vdots  & \vdots  &  \vdots  \\
 44 & \frac{f_8^2 f_{12}^4}{f_4 f_6 f_{24}^2} & 12812 & 2565 & 2565 & 2565 & 2566 & 12805 & 0 & 0 &
   0 & 0 & 7141 & 1219 \\
 \vdots &  \vdots  &  \vdots  &  \vdots  &  \vdots  &  \vdots  &  \vdots  &  \vdots  &  \vdots  &  \vdots  &  \vdots  &  \vdots  & \vdots  &  \vdots  \\
 46 & \frac{f_6^4 f_{12}}{f_3^2 f_{24}} & 12881 & 2603 & 2603 & 2604 & 2603 & 12874 & 1 & 0 & 0 & 0
   & 7141 & 1211 \\
 47 & \frac{f_6 f_8^2 f_{12}}{f_4 f_{24}} & 12812 & 2568 & 2567 & 2567 & 2567 & 12805 & 1 & 0 & 0 &
   0 & 7141 & 1194 \\
 48 & \frac{f_3^2 f_{12}^3}{f_6^2 f_{24}} & 12881 & 2547 & 2547 & 2547 & 2548 & 12874 & 0 & 0 & 0 &
   0 & 7141 & 1184 \\
 49 & \frac{f_4^2 f_{12}^3}{f_6 f_8 f_{24}} & 12812 & 2611 & 2611 & 2611 & 2612 & 12805 & 1 & 0 & 0
   & 0 & 7141 & 1176 \\
  \vdots &  \vdots  &  \vdots  &  \vdots  &  \vdots  &  \vdots  &  \vdots  &  \vdots  &  \vdots  &  \vdots  &  \vdots  &  \vdots  & \vdots  &  \vdots  \\
 58 & \frac{f_4^5 f_{20}^3}{f_2^2 f_8^2 f_{10} f_{40}} & 13070 & 2977 & 2978 & 2977 & 2977 & 13060 &
   8 & 0 & 0 & 0 & 11901 & 2342 \\
   \hline
\end{longtable}
\end{center}
}

\end{document}